\documentclass[11pt]{article}

\usepackage[pagewise]{lineno}
\usepackage[english]{babel}
\usepackage{indentfirst}
\usepackage{latexsym}
\usepackage[colorlinks=true, bookmarksnumbered=true, bookmarksopen=true,
bookmarksopenlevel=3, pdfstartview=FitH, linkcolor=cyan, pdfmenubar=true,
pdftoolbar=true, bookmarks=true,citecolor=cyan, urlcolor=magenta,
filecolor=cyan,menucolor=black,plainpages=false,pdfpagelabels]{hyperref}
\usepackage[lmargin=2cm,tmargin=2cm,bmargin=2cm,rmargin=2cm]{geometry}
\usepackage{amsbsy, amsmath,amsfonts,amssymb,amsthm}
\usepackage{times}
\usepackage{graphicx}
\usepackage{amsmath}
\usepackage{amsfonts}
\usepackage{amssymb}
\usepackage[english]{babel}
\usepackage[latin1]{inputenc}
\usepackage[T1]{fontenc}
\usepackage{amsmath,amssymb,graphics,amsthm}
\usepackage[usenames]{color}

\numberwithin{equation}{section}
\newtheorem{prop}{Proposition}[section]
\newtheorem{defi}[prop]{Definition}
\newtheorem{teo}[prop]{Theorem}

\newtheorem{obs}[prop]{Remark}
\newtheorem{lema}[prop]{Lemma}

\newcommand\supp{\mathop{\rm supp}}
\newcommand\reallywidehat[1]{\savestack{\tmpbox}{\stretchto{  \scaleto{    \scalerel*[\widthof{\ensuremath{#1}}]{\kern-.6pt\bigwedge\kern-.6pt}    {\rule[-\textheight/2]{1ex}{\textheight}}  }{\textheight}}{0.5ex}}\stackon[1pt]{#1}{\tmpbox}}

\newcommand{\vv}{\parallel}

\begin{document}

\title{Global well-posedness for a rescaled Boussinesq system in critical Fourier-Besov spaces}
\author{
{Leithold L. Aurazo-Alvarez$^{1}$}
{\thanks{{Corresponding
author: Leithold L. Aurazo-Alvarez}
\newline
%,lcff@ime.unicamp.br (L.C.F. Ferreira).}\newline
{The author was supported by Cnpq (Brazil) and by University of Campinas, SP, Brazil.
 %L.C.F. Ferreira was supported by FAPESP and CNPq, Brazil.
 }}}\\ 
%EndAName
{\small $^{1}$ IMECC-Department of Mathematics, University of Campinas
}\\ \newline
{\small Campinas, SP, CEP 13083-859, Brazil.}}
\date{}
\maketitle
%%%%%%%%%%%%%%%
%%%%%%%%%%

\begin{abstract}	
In this work we prove a global well-posedness result for a tridimensional rescaled Boussinesq system, with positive full viscosity and diffusivity parameters in the framework of critical Fourier-Besov spaces. This rescaled approach permits to know the behaviour of the system according to relations between both the parameters and the initial velocity and temperature; for instance, it is possible to consider, for small enough viscosity and large diffusivity, a large enough critical Fourier-Besov norm for the initial temperature and it is also possible to consider, for small enough diffusivity and large viscosity, a large enough critical Fourier-Besov norm for both the velocity and the temperature.
\end{abstract}

\renewcommand{\abstractname}{Abstract}

\section{Introduction}

In this work we deal with a particular case of the following family of Boussinesq equations in $N$ dimensions,

\begin{equation*}
(BQ)^{(\nu_1,\nu_2\dots, \nu_N)}_{(\eta_1,\eta_2,\dots,\eta_N)}\,\,\,\left\{
\begin{split}\label{BQ}
&  \partial_{t}u + (u\cdot\nabla)u - \nu_1\partial^{2}_{x_1}u - \nu_2\partial^{2}_{x_2}u-\dots - \nu_N\partial^{2}_{x_N}u  
+ \nabla P = \theta e_N,\\
&  \partial_{t}\theta + (u\cdot\nabla)\theta - k_1\partial^{2}_{x_1}\theta 
- k_2\partial^{2}_{x_2}\theta - \dots - k_N\partial^{2}_{x_N}\theta= 0,\\
&  \mbox{div}\,u=0,\,\,\mbox{for}\,\,x=(x_1, x_2, \dots, x_N)\in\mathbb{R}^{N}\,\,\mbox{and}\,\, t>0,\,\,\mbox{and} \\
&  u(x,0)=u_{0}(x),\,\,\theta(x,0)=\theta_{0}(x),\,\,\mbox{for}\,x\in
\mathbb{R}^{N},
\end{split}
\right.  
\end{equation*}
%%%%%%%%%%%%%%%%%%%%%%%%%%%%%%%%%%%%%%%%%
where $u=u(x,t)=(u^{1}(x,t), u^{2}(x,t),\dots, u^{N}(x,t))$ represents the $N$-dimensional
velocity vector field, the scalar quantity
$\theta=\theta(x,t)$ denotes, depending of the
situation, for instance, if we deal with thermal
convection, the temperature vatiation in a
gravity field, the scalar function $P=P(x,t)$
denotes the pressure, each one of the
parameters, $\nu_i$ and $k_i$, represents the
ith - viscosity and the ith - thermal
diffusivity, respectively,
for $i=1, 2, \dots, N$, and the 
term $\theta e_N$
denotes the buoyancy force. Boussinesq equations describes the 
interaction of thermal convection with
the dynamics of the fluid, serve as a fundamental model in geophysical
phenomena as ocean circulation, atmospheric fronts and situations where rotation
and stratification play a relevant role
(see \cite{PedloskyGeophysicalFD}, \cite{BatchelorAnIntrodFD} and \cite{MajdaIntrodPDEWAO}).

In the rest of this section we recall some previous studies on this system and related ones. We start with two-dimensional models. In the work (\cite{GuoSpectralmethodNBsystem}, 1989) G. Boling 
considered the following
periodical Newton-Boussinesq equations

\begin{equation}\label{NBQ}
(NB)\,\,\,\left\{
\begin{split}
&  \partial_{t}\zeta+u\partial_x\zeta+v\partial_y\zeta-\Delta \zeta
= -\frac{R_a}{p_r}\partial_x\theta,\\
&  \Delta\psi=\zeta,\,\,u=\psi_y,\,\,v=-\psi_x,\\
&\partial_t\theta+u\partial_x\theta+v\partial_y\theta=\frac{1}{p_r}\Delta\theta,
\end{split}
\right. 
\end{equation}
where the vector $(u,v)$ represents the
velocity, the scalar $\theta$ represents
the temperature, the term $\psi$ is the
flow function, the scalar $\zeta$ is the vortex,
the parameter $p_r$ is the Prandtl number
and the parameter $R_a$ is the Rayleigh number. The author rewrites the system (\ref{NBQ})
as follow

\begin{equation}\label{NBQequivalente}
(NB^{*})\,\,\,\left\{
\begin{split}
&  \partial_{t}\Delta\psi +J(\psi, \Delta \psi)-\Delta^{2}\psi=-\frac{R_a}{p_r}\frac{\partial\theta}{\partial x}\\
&\partial_{t}\theta + J(\psi,\theta)=\frac{1}{p_r}\Delta\theta,
\end{split}
\right.
\end{equation}
%%%%%%%%%%%%%%%%
where $J(u,v)=u_yv_x-u_xv_y$ and he is looking for mean zero $2\pi$-periodical solutions $\psi(x,y,t)$ and $\theta(x,y,t)$, considering periodical initial data $\psi_0$ and $\theta_0$ defined in $\Omega=(0, 2\pi)^{2}$. He studied this system (\ref{NBQequivalente}) in order to get the existence and uniqueness of global generalized solution as well as its uniqueness. In particular, he proved the existence of generalized solutions $\psi$ and $\theta$ in the respective spaces $\psi\in L^{\infty}(0,T; H^{2}(\Omega))$ and $\theta\in L^{\infty}(0,T; H^{1}(\Omega))$. The assumption used to get these results are $\psi_0\in H^{1}(\Omega)$ and $\theta_0\in L^{2}(\Omega)$ (or $\psi_0\in H^{2}(\Omega)$ and $\theta_0\in H^{1}(\Omega)$ ) and an inequality involving the parameter $p_r$. In the work (\cite{Hou-LiGlobalWVBQ}, 2005) T.Y. Hou and C. Li proved the existence of a unique smooth global solution $(u,\theta)$ for the system $(BQ)^{(\nu,\nu)}_{(0, 0)}$, with $\nu=1$, assuming arbitrarily large initial data in Sobolev spaces, $u_0\in H^{m}(\mathbb{R}^{2})$ and $\theta_0\in H^{m-1}(\mathbb{R}^{2})$, for $m\geq 3$. In the work (\cite{ChaeGlobalRegBQPartViscosity}, 2006)
 D. Chae proved the existence of a unique
 solution $(u,\theta)$ for the system
 $(BQ)^{(\nu, \nu)}_{(0, 0)}$, with $\nu>0$, in the space
 
 \begin{equation*}
\begin{split}
&u\in C([0,\infty); H^{m}(\mathbb{R}^{2}))
\cap C(0,T; H^{m+1}(\mathbb{R}^{2}))\,\,\mbox{and}\\
&\theta\in C([0,\infty); H^{m}(\mathbb{R}^{2})),
\end{split}     
 \end{equation*}
 considering arbitrarily initial data $(u_0, \theta_0)\in H^{m}(\mathbb{R}^{2})$, for $m>2$. He established the convergence of solutions for the Boussinesq system $(BQ)^{(\nu,\nu)}_{(k, k)}$, with $\nu>0$ and $k>0$, to the corresponding solutions of the systems $(BQ)^{(\nu,\nu)}_{(0,0)}$ in the space $C([0,T]; H^{s}(\mathbb{R}^{2}))$, for $s<m$, as $k\rightarrow 0$. He also proved the existence of a unique solution $(u,\theta)$ for the system $(BQ)^{(0,0)}_{(k,k)}$, with $k>0$, in the space 
 
\begin{equation*}
\begin{split}
&u\in C([0,\infty); H^{m}(\mathbb{R}^{2}))\,\,\mbox{and}\\
&\theta\in C([0,\infty); H^{m}(\mathbb{R}^{2}))\cap L^{2}(0,T; H^{m+1}(\mathbb{R}^{2})),
\end{split}     
 \end{equation*}
 considering initial data $(u_0,\theta_0)\in H^{m}(\mathbb{R}^{2})$, for $m>2$. Moreover, he established the convergence of solutions of the systems 
 $(BQ)^{(\nu,\nu)}_{(k,k)}$, with $\nu>0$ and $k>0$, to the corresponding solutions of the system $(BQ)^{(0,0)}_{(k,k)}$, in the space $C([0,T]; H^{s}(\mathbb{R}^{2}))$, for $s<m$, as $\nu\rightarrow 0$. In the work (\cite{HmidiKeraaniOntheGlobalBQZeroDiffus}, 2007) T. Hmidi and S. Keraani proved the existence of a global weak solution $(u,\theta)$ for the system $(BQ)^{(\nu,\nu)}_{(0,0)}$, with $\nu>0$, in the space

\begin{equation*}
\begin{split}
&u\in C(\mathbb{R}_{+}; H^{s}(\mathbb{R}^{2}))\cap L^{2}_{loc}(\mathbb{R}_{+}; H^{\min(s+1,2)}(\mathbb{R}^{2}))\,\,\mbox{and}\\
&\theta\in C(\mathbb{R}_{+}; L^{2}(\mathbb{R}^{2})),
\end{split}     
 \end{equation*}
considering initial data $u_0\in H^{s}(\mathbb{R}^{2})$, with $0\leq s<2$, and $\theta_0\in L^{2}(\mathbb{R}^{2})$. They also proved the uniqueness of these weak solutions by assuming additional regularity for the initial data, for instance, $u_{0}\in H^{s}(\mathbb{R}^{2})$ and $\theta_{0}\in B^{0}_{2,1}\cap B^{0}_{p,\infty}(\mathbb{R}^{2})$, for $0<s\leq 2$ and $2<p\leq \infty$. In the work (\cite{AbidiHmidiOntheGlobalBQ}, 2007) H. Abidi and T. Hmidi proved the existence of a unique global solutions $(u, \theta, P)$ for the system $(BQ)^{(\nu,\nu)}_{(0, 0)}$, with $\nu>0$, such that

\begin{equation*}
\begin{split}
&u\in C(\mathbb{R}_{+}; L^{2}\cap B^{-1}_{\infty, 1}(\mathbb{R}^{2}))\cap L^{2}_{loc}(\mathbb{R}_{+}; H^{1}(\mathbb{R}^{2}))\cap L^{1}_{loc}(\mathbb{R}_{+}; B^{1}_{\infty, 1}(\mathbb{R}^{2}))\\
&\theta\in C(\mathbb{R}_{+}; B^{0}_{2,1}(\mathbb{R}^{2}))\,\,\mbox{and}\,\,
\nabla P\in L^{1}_{loc}(\mathbb{R}_{+}; B^{0}_{2,1}(\mathbb{R}^{2})),
\end{split}     
 \end{equation*}
taking initial data $u_0\in L^{2}\cap B^{-1}_{\infty,1}(\mathbb{R}^{2})$ and 
$\theta_0\in B^{0}_{2,1}(\mathbb{R}^{2})$. In the work (\cite{DanchinPaicuGlobalIBQYudovichData}, 2009) R. Danchin and 
M. Paicu proved the existence of a unique
global solution $(u,\theta)$ for the system
$(BQ)^{(0,0)}_{(k,k)}$, with $k>0$, in the space

\begin{equation*}
\begin{split}
&u\in C^{0,1}_{loc}(\mathbb{R}_{+}; L^{2}(\mathbb{R}^{2})),\,\,\omega\in L^{\infty}_{loc}(\mathbb{R}_{+}; L^{r}\cap L^{\infty}(\mathbb{R}^{2}))\,\,\mbox{and}\\
&\theta \in C(\mathbb{R}_{+}; L^{2}\cap B^{-1}_{\infty, 1}(\mathbb{R}^{2}))\cap L^{2}_{loc}(\mathbb{R}_{+}; H^{1}(\mathbb{R}^{2}))\cap L^{1}_{loc}(\mathbb{R}_{+}; B^{1}_{\infty, 1}(\mathbb{R}^{2})),
\end{split}     
 \end{equation*}
where $\omega$ denotes the vorticity and they are considering initial data 
$u_0\in L^{2}(\mathbb{R}^{2})$, $\omega_0\in L^{r}\cap L^{\infty}(\mathbb{R}^{2})$, for some $2\leq r <\infty$, and $\theta_0\in L^{2}\cap B^{-1}_{\infty, 1}(\mathbb{R}^{2})$. In the work (\cite{DanchinPaicuGlobalExAnisBQ}, 2011) R. Danchin and M. Paicu proved the existence of a unique global solution $(u,\theta)$ for the system 
$(BQ)^{(\nu_1,0)}_{(0,0)}$, with $\nu_1>0$, such that $u\in C_{w}(\mathbb{R}_{+}; H^{1}(\mathbb{R}^{2}))$, $\omega\in L^{\infty}_{loc}(\mathbb{R}_{+};\sqrt{L})$ and $\theta\in C_{w}(\mathbb{R}_{+}; L^{\infty})\cap C(\mathbb{R}_{+}; H^{s-\varepsilon}(\mathbb{R}^{2}))$, for all $\varepsilon>0$, and considering initial data $u_0\in H^{1}(\mathbb{R}^{2})$, $\omega_0\in\sqrt{L}$ and $\theta_0\in H^{s}\cap L^{\infty}(\mathbb{R}^{2})$, for $\frac{1}{2}<s\leq 1$, where 

\begin{equation*}
\parallel f\parallel_{\sqrt{L}}=\displaystyle{\sup_{p\geq 2}}\frac{\parallel f\parallel_{L^{p}(\mathbb{R}^{2})}}{\sqrt{p-1}}.    
\end{equation*}
They also proved the existence of a unique solution $(u,\theta)\in C_{w}(\mathbb{R}_{+}; H^{1}(\mathbb{R}^{2}))$ for the system $(BQ)^{(0,0)}_{(k_1,0)}$, with $k_1>0$, and such that  
$\omega\in L^{\infty}_{loc}(\mathbb{R}_{+}; L^{\infty}(\mathbb{R}^{2}))$, 
$\theta\in L^{\infty}(\mathbb{R}_{+}; H^{1}(\mathbb{R}^{2}))$ and $\partial_1\theta\in L^{2}(\mathbb{R}_{+}; H^{1}(\mathbb{R}^{2}))$,
where the initial data they considered have the next regularities, $u_0\in H^{1}(\mathbb{R}^{2})$, $\omega_0\in L^{\infty}(\mathbb{R}^{2})$, $\theta_0\in H^{1}(\mathbb{R}^{2})$ and $\mid\partial_1\mid^{s}\theta_0\in L^{2}(\mathbb{R}^{2})$, for $1<s<\frac{3}{2}$, where $(\mid\partial_1\mid^{s}\theta_0)^{\wedge}(\xi)=\mid \xi_1\mid^{s}\hat{\theta}_0(\xi)$. In the work (\cite{CaoWuGlobalRegVerticalDiss}, 2013) C. Cao and J. Wu proved for initial data $(u_0, \theta_0)\in H^{2}(\mathbb{R}^{2})$ and for each $T>0$, the existence of a unique classical solution $(u,\theta)\in C([0,T]; H^{2}(\mathbb{R}^{2}))$ for the Boussinesq system $(BQ)^{(0,\nu_2)}_{(0,k_2)}$, with $\nu_2>0$ and $k_2>0$. In the work (\cite{LiTitiGlobalBQVerticalDiss}, 2016) J. Li and E.S. Titi proved, for any initial data $(u_0,\theta_0)$ satisfying $(u_0, \theta_0, \partial_{x}u_0,\partial_{x}\theta_0)\in L^{2}(\mathbb{R}^{2})$ and any value for $T\in (0,\infty)$, the existence of a solution $(u,\theta)\in 
C([0,T]; L^{2}(\mathbb{R}^{2}))\cap L^{2}(0,T; H^{1}(\mathbb{R}^{2}))$ for the Boussinesq system $(BQ)^{(0,\nu_2)}_{(0,k_2)}$, with $\nu_2>0$ and $k_2>0$, as well as the uniqueness of this solution if it belongs to some spaces with determined regularities.

%%%%%%%%%%%%%%%%%%%%%%%%%
%%%%%%%%%%%%%%%%%%%%%%%%%%
%%%%%%%%%%%%%%%%%%%%%%%%%%%
Now, we recall some previous results on these Boussinesq equations in general $N$-dimensional spaces, mostly for $N\geq 3$. In the work (\cite{Cannon-DiBenedetto}, 1980) J.R. Cannon and E. DiBenedetto considered the following Boussinesq equations in the infinite cylinder $S_T=\mathbb{R}^{N}\times (0,T]$

\begin{equation}\label{CBQ}
(CBQ)\,\,\,\left\{
\begin{split}
&  \frac{\partial\theta}{\partial t}-\Delta\theta+(u\cdot \nabla)\theta =q(x,t)\\
& \theta(x,0)=\theta_0(x)\\
&\frac{\partial u}{\partial t}-\Delta u+(u\cdot\nabla) u+\nabla P=f(\theta(x,t))\\
&\nabla\cdot u=0\\
&u(x,0)=u_0(x),
\end{split}
\right.
\end{equation}
where $x\in\mathbb{R}^{N}$ and $t\in (0,T]$, for a fixed time $T>0$, and $
N$ a positive entire number. Here $q(\cdot)$ and $f(\cdot)$ are known data functions and $P(\cdot)$ is the pressure term. Also, the term $u=(u_1(x,t),\dots,u_N(x,t))$ is the velocity  vector field and $\theta=\theta(x,t)$ represents the scalar temperature. They are interested in looking for solutions $U=(u,\theta)$ in the space

\begin{equation*}
L^{p,q}_{N+1}(S_T)=\lbrace
U:S_T\longrightarrow\mathbb{R}^{N+1}; \parallel U\parallel_{p,q}=
\displaystyle\sum_{j=1}^{N+1}\lbrace 
\int_{0}^{T}\left[ \int_{\mathbb{R}^{N}}\mid U_j\mid^{p}\,dx\right]^{\frac{q}{p}}\,dt
\rbrace^{\frac{1}{q}}<\infty
\rbrace,    
\end{equation*}
for $p,q\geq 2$, $\frac{N}{p}+\frac{2}{q}\leq 1$ and $N<p<\infty$. Here 
$(U_1,\dots , U_N)=u$ and $U_{N+1}=\theta$. 
Firstly, they obtain a non-linear integral equation whose solutions in these spaces $L^{p,q}_{N+1}(S_T)$ are equivalent to the weak solutions of the Boussinesq system (CBQ) in the same kind of spaces. This integral equation can be expressed as follows

\begin{equation}\label{CBQIntegral}
U(x,t)+B(U,U)=U_0(x,t)+F(U)(x,t),    
\end{equation}
where $U(x,t)=(u(x,t),\theta(x,t))$, the term $U_0$ depends on the given initial velocity $u_0$ and temperature $\theta_0$, the term $F(U)(x,t)$ depends on the functions $f(\theta(y,s))$ and $q(y,s)$ and $B(U,V)(x,t)$ is the bilinear term depending of $U=(u,\theta)$ and $V=(\tilde{u},\tilde{\theta})$. They proved the local existence of a unique solution $U\in L^{p,q}_{N+1}(S_T)$ for the system (\ref{CBQIntegral}), assuming $\frac{N}{p}+\frac{2}{q}\leq 1$, $N<p<\infty$ and an initial data $U_0\in L^{r}_{N+1}(\mathbb{R}^{N})$, for $0<\frac{N}{r}<\frac{2}{q}+\frac{N}{p}\leq 1$. Moreover, if consider small initial data $U_0\in L^{r_1}_{N+1}(\mathbb{R}^{N})\cap L^{r_2}_{N+1}(\mathbb{R}^{N})$, for suitable $r_1$ and $r_2$ related to $p,q$ and $N$, they obtained, for each $T>0$, a unique solution $U\in L^{p,q}_{N+1}(S_T)$ of the nonlinear integral equation (\ref{CBQIntegral}). They also obtained a result of permanency of regularity, that is, if the initial data has the regularity $\frac{\partial}{\partial x_k}U_0\in L^{p/2,q/2}_{N+1}(S_T)$ then the solution also has the same regularity
$\frac{\partial}{\partial x_k}U\in L^{p/2,q/2}_{N+1}(S_T)$. In the work (\cite{DanchinPaicutheoremesLFK}, 2008) R. Danchin and M. Paicu considered the $N$-dimensional Boussinesq system with partial viscosity $(BQ)^{(\nu_1, \dots, \nu_N)}_{(k_1, \dots, k_N)}$, with $\nu_1=\dots=\nu_N=\nu>0$ and $k_1=\dots=k_N=k=0$. Firstly they prove the existence of a global weak solution $(\theta, u)$ for the system $(BQ)^{(\nu_1, \dots, \nu_N)}_{(k_1, \dots, k_N)}$, such that $\theta \in L^{\infty}(\mathbb{R}_{+};L^{p}(\mathbb{R}^{N}))$ and $u\in L^{\infty}_{loc}(\mathbb{R}_{+};L^{2}(\mathbb{R}^{N}))\cap L^{2}_{loc}(\mathbb{R}_{+};H^{1}(\mathbb{R}^{N}))$, whenever the initial data are in the next spaces, $\theta_0\in L^{p}(\mathbb{R}^{N})$, with 
$\frac{2N}{N+2}<p\leq 2$, and $u_0\in L^{2}(\mathbb{R}^{N})$. They also proved the existence of a unique local solution $(\theta, u, \nabla P)$ in the space

\begin{equation*}
E_T={\cal C}([0,T]; \dot{B}^{0}_{N,1})\times ( {\cal C}([0,T]; \dot{B}^{0}_{N,1})\cap L^{1}(0,T; \dot{B}^{2}_{N,1}))\times L^{1}(0,T;\dot{B}^{0}_{N,1}),     
\end{equation*}
considering initial data $\theta_0\in \dot{B}^{0}_{N,1}$ and $u_0\in \dot{B}^{0}_{N,1}$, for $N\geq 2$. Moreover, they established a continuation criteria for this local solution. They also proved the existence of a unique global solution $(\theta, u,\nabla P)$ satisfying 

\begin{equation*}
\begin{split}
&(\theta,u)\in {\cal C}(\mathbb{R}_{+};\dot{B}^{0}_{N,1})\times
{\cal C}(\mathbb{R}_{+};\dot{B}^{\frac{N}{p}-1}_{p,1})\cap L^{1}_{loc}(\mathbb{R}_{+};\dot{B}^{\frac{N}{p}+1}_{p,1})\,\,\mbox{and}\\
&\nabla P\in L^{1}_{loc}(\mathbb{R}_{+};\dot{B}^{\frac{N}{p}-1}_{p,1}),
\end{split}
\end{equation*}
whenever the initial data $\theta_0\in\dot{B}^{0}_{N,1}$ and $u_0\in \dot{B}^{1+\frac{N}{p}}_{p,1}$, for $N\leq p\leq \infty$, satisfy a smallness condition 

\begin{equation*}
\parallel u_0\parallel_{L^{N,\infty}}+\nu^{-1}\parallel \theta_0\parallel_{L^{\frac{N}{3}}}\leq c\nu,    
\end{equation*}
for some constant $c>0$ independent of $N$. Here $L^{N,\infty}$ denotes a Lorentz space. Finally, they proved the existence of a unique global solutions $(\theta, u,\nabla P)$  satisfying

\begin{equation*}
\theta\in {\cal C}(\mathbb{R}_{+};L^{2})\,\,\mbox{and}\,\,u\in {\cal C}(\mathbb{R}_{+};L^{2})\cap L^{2}_{loc}(\mathbb{R}_{+};H^{1}),
\end{equation*}
for any initial data $(\theta_0,u_0)\in L^{2}(\mathbb{R}^{2})$. In the work (\cite{DanchinPaicuExistUniqueLorentz}, 2008) R. Danchin and M. Paicu studied the $N$-dimensional Boussinesq system with partial viscosity $(BQ)^{(\nu_1, \dots, \nu_N)}_{(k_1, \dots, k_N)}$, with $\nu_1=\dots=\nu_N=\nu>0$ and $k_1=\dots=k_N=k=0$. They proved the existence of a global solution $(\theta, u,\nabla P)$ for this system,  satisfying

\begin{equation*}
\begin{split}
&u\in L^{\infty}(\mathbb{R}_{+};L^{N,\infty})\,\,\mbox{and}\\
&\theta\in L^{\infty}(\mathbb{R}_{+};L^{1}\cap L^{p,\infty})\,\,\mbox{if}\,\,N=3, \,\,\theta\in L^{\infty}(\mathbb{R}_{+};L^{\frac{N}{3},\infty}\cap L^{p,\infty})\,\,\mbox{if}\,\,N\geq 4,
\end{split}    
\end{equation*}
for small initial condition

\begin{equation*}
\parallel u_0\parallel_{L^{3,\infty}}+\parallel\theta_0\parallel_{L^{1}}\leq c\nu\,\,\mbox{if}\,\,N=3, \,\, \parallel u_0\parallel_{L^{N,\infty}}+\parallel\theta_0\parallel_{L^{\frac{N}{3},\infty}}\leq c\nu\,\,\mbox{if}\,\,N\geq 4,    
\end{equation*}
where $c>0$ is a constant depending only on $N$, and under the additional assumption $u_0\in L^{p,\infty}$, for large enough value of $p$. They established a uniqueness result for the Boussinesq system $(BQ)^{(\nu_1, \dots, \nu_N)}_{(k_1, \dots, k_N)}$, for dimension $N\geq 3$, assuming that solutions $(\theta_1,u_1,\nabla P_1)$ and $(\theta_2,u_2,\nabla P_2)$ are in the following spaces

\begin{equation*}
\theta_{i}\in L^{\infty}_{T}(\dot{B}^{-1+\frac{N}{P}}_{p,\infty})\,\,\mbox{and}\,\,u_i\in L^{\infty}_{T}(\dot{B}^{-1+\frac{N}{P}}_{p,\infty})\cap 
\tilde{L}^{1}_{T}(\dot{B}^{1+\frac{N}{P}}_{p,\infty}),\,\,\mbox{for}\,i=1,2,
\end{equation*}
and it is valid the condition

\begin{equation*}
\parallel u_1\parallel_{\tilde{L}^{1}_{T}(\dot{B}^{1+\frac{N}{P}}_{p,\infty})}+\nu^{-1}\parallel u_2\parallel_{L^{\infty}_{T}(\dot{B}^{-1+\frac{N}{P}}_{p,\infty})} \leq c,  
\end{equation*}
for some constant $c>0$ depending only on $N$ and $p$. Finally, for the limit case, $p=2N$, they also proved the uniqueness of a global solution $(\theta,u,\nabla P)$ with 

\begin{equation*}
\begin{split}
&u\in L^{\infty}(\mathbb{R}_{+};L^{N,\infty})\,\,\mbox{and}\\  
&\theta\in L^{\infty}(\mathbb{R}_{+};L^{\frac{N}{3},\infty})\,\,\mbox{if}\,\,N\geq 4,\,\,\theta\in L^{1}(\mathbb{R}_{+};L^{1})\,\,\mbox{if}\,\,N=3.
\end{split}    
\end{equation*}
In the work (\cite{Ferreira-Villamizar-R-Existence-convection-pseudomeasure-2008}, 2008) L.C.F. Ferreira and E.J. Villamizar-Roa studied a $N$-dimensional generalized convection problem with gravitacional field depending of the space variable, in the context of ${\cal PM}^{a}$ - spaces. They obtained local and global well-posedness result for this system and, in particular, for the  B\'enard problem, which correspond to a generalized Newtonian gravitational field. They also study regularizations properties for solutions. In the work (\cite{HmidiRoussetGlobalwellNSBQaxi}, 2010) T. Hmidi and F. Rousset considered the Navier-Stokes-Boussinesq system 
$(BQ)^{(\nu, \nu, \nu)}_{(k, k ,k)}$, with $\nu=1$ and $k\geq 0$. They established the existence of a unique global solutions $(u,\theta)$ for the above system $(BQ)^{(\nu, \nu, \nu)}_{(k, k ,k)}$, satisfying

\begin{equation*}
\begin{split}
&u\in{\cal C}(\mathbb{R}_{+}; H^{1})\cap L^{2}_{loc}(\mathbb{R}_{+};H^{2})\cap 
L^{1}_{loc}(\mathbb{R}_{+};B^{1}_{\infty,1}),\\
&\theta\in{\cal C}(\mathbb{R}_{+};L^{2}\cap B^{0}_{3,1})\,\,\mbox{and}\,\,\frac{\omega}{r}\in L^{\infty}_{loc}(\mathbb{R}_{+}),
\end{split}    
\end{equation*}
 assuming $u_0\in H^{1}$ a divergence free axisymmetric vector field, $\frac{\omega_0}{r}\in L^{2}$ and $\theta_0\in L^{2}\cap L^{m}$, for $m>3$, (or $\theta_0\in L^{2}\cap \dot{B}^{0}_{3,1}$) an axisymmetric function.  Here there is no smallness assumption on the initial data. In the work (\cite{Ferreira-Villamizar-R-Stability-Boussinesq-Weak-Lp-2010}, 2010) L.C.F. Ferreira and E.J. Villamizar-Roa studied the problem of stability of steady solutions of the Boussinesq system, with positive parameters, in the context of weak-$L^{p}$ spaces, for several types of domains, as the entire space $\mathbb{R}^{N}$, the half space $\mathbb{R}^{N}_{+}$, bounded or exterior domains with some smoothness for the boundary. Thus, by considering small initial perturbations in the weak-$L^{N}$ norm, a control in time for the weak-$L^{b}$ norm of the external force, for some $b$,  and decaying properties for the initial perturbations involving the associated semigroup, they proved that non-steady solutions for the perturbed system converges to the steady ones in some weak-$L^{q}$ spaces, for $3\leq N<q<\infty$. In the work (\cite{AbidiHmidiKeraaniGlobalregNSB}, 2011) H. Abidi, T. Hmidi and S. Keraani studied the tridimensional Boussinesq system 
$(BQ)^{(\nu, \nu, \nu)}_{(0, 0, 0)}$, with $\nu=1$, for axisymmetric initial data. One important quantity in order to establish the regularity of the Boussinesq system and related systems is the control in an appropiated norm of the called vorticity, which in cylindrical coordinates is given by the expression $\omega=(\partial_{z} u^{r}-\partial_{r}u^{z})e_{\theta}=\omega_{\theta}e_{\theta}$.
The authors proved the existence of a unique global solution $(u,\theta)$ such that

\begin{equation*}
\begin{split}
&u\in C(\mathbb{R}_{+};H^{1})\cap L^{1}_{loc}(\mathbb{R}_{+};W^{1,\infty}),\\
&\theta\in L^{\infty}_{loc}(\mathbb{R}_{+};L^{2}\cap L^{\infty})\,\,\mbox{and}\,\,
\frac{\omega}{r}\in L^{\infty}_{loc}(\mathbb{R}_{+};L^{2}),
\end{split}    
\end{equation*}
assuming a divergence free initial axisymmetric velocity $u_0\in H^{1}$, the quantity $\frac{\omega_{0}}{r}\in L^{2}$, the temperature $\theta_0\in L^{2}\cap L^{\infty}$ depending only on $(r,z)$ and the support of $\theta_0$ disjoint of the axis $(OZ)$ with their projection on this axis being compact. In the work (\cite{Ferreira-Almeida-Wellposedness-Boussinesq-Morrey-2011}, 2011) L.C.F. Ferreira and M.F. de Almeida studied the existence and the asymptotic behaviour of global solutions in Morrey spaces ${\cal M}_{p,\lambda}$ for the Boussinesq system, with positive parameters and external force. They proved the global well-posedness for the system, with small initial data in ${\cal M}_{p,\lambda}$ - spaces and a small enough external force in time-dependent spaces based on Morrey spaces. In particular, this space for the external force permits to consider the Newtonian gravitational field and thus cover the B\'enard problem. They also study the existence of self-similar solutions. In the work (\cite{MiaoZhengBQHorizontalDissipation}, 2013) C. Miao and X. Zheng studied the following Boussinesq system with horizontal viscosity and horizontal diffusion $(BQ)^{(\nu, \nu, 0)}_{(k, k , 0)}$, with $\nu=1$ and $k>0$. They proved the existence of a unique global solution $(u,\theta)$ for the system $(BQ)^{(\nu, \nu, 0)}_{(k, k ,0)}$, in the spaces
 
 \begin{equation*}
\begin{split}
&u\in{\cal C}(\mathbb{R}_{+};H^{1})\cap L^{2}_{loc}(\mathbb{R}_{+};H^{1,2}\cap H^{2,1}),\\
&\theta\in {\cal C}(\mathbb{R}_{+};H^{0,1})\cap L^{2}_{loc}(\mathbb{R}_{+};H^{1,1})\,\,\mbox{and}\\
&\frac{\omega}{r}\in L^{\infty}_{loc}(\mathbb{R}_{+},L^{2})\cap L^{2}_{loc}(\mathbb{R}_{+}, H^{1,0}),
\end{split}     
 \end{equation*}
 where they considered a divergence free axisymmetric initial vector field $u_0\in H^{1}$, $\frac{\omega_0}{r}\in L^{2}$, $\partial_{z}\omega_0\in L^{2}$ and $\theta_0\in H^{0,1}$ an axisymmetric scalar function. They also obtained the existence of a unique global solution $(\theta, u)$ for the system $(BQ)^{(\nu, \nu, 0)}_{(k, k ,0)}$, in the space
 
  \begin{equation*}
\begin{split}
&u\in{\cal C}_{w}(\mathbb{R}_{+};H^{1})\cap L^{2}_{loc}(\mathbb{R}_{+};H^{1,2}\cap H^{2,1}),\\
&\theta\in {\cal C}_{w}(\mathbb{R}_{+};H^{0,1})\cap L^{2}_{loc}(\mathbb{R}_{+};H^{1,1})\cap {\cal C}(\mathbb{R}_{+};L^{2})\,\,\mbox{and}\\
&\frac{\omega}{r}\in L^{\infty}_{loc}(\mathbb{R}_{+}; L^{2})\cap L^{2}_{loc}(\mathbb{R}_{+}; H^{1,0}),
\end{split}     
 \end{equation*}
assuming more general initial data, in the sense that the initial vorticity $\omega_0\in L^{\infty}$ and all the another initial conditions remains in the same spaces. In the work (\cite{MiaoZhengBQAxihorizontalViscosity}, 2014) C. Miao and X. Zheng considered the Boussinesq sytem with horizontal viscosity $(BQ)^{(\nu, \nu, 0)}_{(0, 0 ,0)}$, with $\nu>0$. They proved the existence of a unique global solution $(\theta, u)$ such that 
 
 \begin{equation*}
\begin{split}
&u\in C(\mathbb{R}_{+};H^{1}(\mathbb{R}^{3}))\,\,\mbox{and}\\
&\theta\in L^{\infty}(\mathbb{R}^{+}; L^{\infty}(\mathbb{R}^{3}))\cap 
C(\mathbb{R}^{+}; H^{1}(\mathbb{R}^{3})),
\end{split}     
 \end{equation*}
 by assuming the initial divergence free axisymmetric vector field $u_0\in H^{1}(\mathbb{R}^{3})$, the term $\frac{\omega_{0}}{r}\in L^{2}(\mathbb{R}^{3})$, $\partial_{z}\omega_0\in L^{2}$ and the initial axisymmetric scalar function $\theta_0\in H^{1}(\mathbb{R}^{3})\cap L^{\infty}(\mathbb{R}^{3})$, with both the conditions $\supp{\theta_0}\cap (OZ) =\phi$ and $\Pi_{z}(\supp{\theta_0})$ being compact. In the work (\cite{FangLeZhangGlobalAxiNonzeroSwirl}, 2018) D. Fang, W. Le and T. Zhang studied the Boussinesq system $(BQ)^{(\nu, \nu, \nu)}_{(k, k ,k)}$, with $\nu>0$ and $k\geq0$. The authors proved the local existence of a unique solution $(u,\theta)$ 
for the Boussinesq system in the space

\begin{equation*}
\begin{split}
&u\in C([0,T]; H^{2}(\mathbb{R}^{3}))\cap \tilde{L}^{1}([0,T]; \dot{H}^{4}(\mathbb{R}^{3}))\,\,\mbox{and}\\
&\theta\in C([0,T]; L^{2}(\mathbb{R}^{3})\cap B^{0}_{3,1}(\mathbb{R}^{3}))
\cap \tilde{L}^{1}([0,T];B^{2}_{2,2}\cap B^{2}_{3,1}),\,\,\mbox{for}\,\,k>0,
\end{split}
\end{equation*}
and in the space

\begin{equation*}
\begin{split}
&u\in C([0,T]; H^{2}(\mathbb{R}^{3}))\cap L^{2}([0,T]; \dot{H}^{3}(\mathbb{R}^{3}))\,\,\mbox{and}\\
&\theta\in C([0,T]; H^{1}(\mathbb{R}^{3})),\,\,\mbox{for}\,\,k=0,
\end{split}
\end{equation*}
under the initial condition $(u_0, \theta_0)$ satisfying one of the assumption
\begin{itemize}
\item [i.]    
$
u_0\in H^{2}(\mathbb{R}^{3}),\,\,\theta_0\in L^{2}(\mathbb{R}^{3})\cap 
B^{0}_{3,1}(\mathbb{R}^{3}),\,\,\mbox{when}\,\,k>0,
$
or 
\item [ii.]
$
u_0\in H^{2}(\mathbb{R}^{3}),\,\,\theta_0\in H^{1}(\mathbb{R}^{3})\cap 
B^{0}_{3,1}(\mathbb{R}^{3})\,,\,\mbox{when}\,\,k=0.
$
\end{itemize}
They also established the existence of a unique global solution $(u,\theta)$ for the Boussinesq system under some additional conditions, for instance $ru^{\theta}_0\in L^{1}\cap L^{4}(\mathbb{R}^{3})$ and a smallness condition for the norm $\parallel r^{\varepsilon}u^{\theta}_0\parallel_{L^{\frac{3}{1-d}}}$, for $0\leq \varepsilon<1$. Finally, they established a decay estimative for this global solution, in particular, they proved that
 
\begin{equation*}
\parallel \theta(t)\parallel^{2}_{L^{2}(\mathbb{R}^{3})}\lesssim\, <t>^{-\frac{3}{4}},\,\,\parallel u(t)\parallel^{2}_{L^{2}(\mathbb{R}^{3})}\lesssim\, <t>^{\frac{1}{2}}\,\,\mbox{and}\,\,\parallel ru^{\theta}(t)\parallel^{2}_{L^{2}(\mathbb{R}^{3})}\lesssim\, <t>^{-\frac{3}{2}},    
\end{equation*}
 where the symbol $<t>=\sqrt{1+t}$. In the work (\cite{Ferreira-Perez-Villamizar-B-Lorentz-M-SBQ-2018}, 2018) L.C.F. Ferreira, J.E. P\'erez-L\'opez and E.J. Villamizar-Roa studied the $N$-dimensional stationary Boussinesq-Coriolis with external forces, in the context of Besov-Lorentz-Morrey and Besov-weak-Morrey spaces. They proved the existence of stationary solutions for small external forces in these spaces. In particular, this spaces allow to consider the gravitacional field corresponding to the B\'enard problem. In the work (\cite{Aurazo-Ferreira-Global-BC-Stratification-FBM-2021}, 2021) L.C.F. Ferreira and the author studied the tridimensional fractional Boussinesq-Coriolis system with stratification in the context of critical Fourier-Besov-Morrey spaces. We proved a global well-posedness result for small
initial data in this spaces, for a range of values for the exponent of minus the Laplacian. There, we also consider
the critical case for this system and for the tridimensional fractional Navier-Stokes-Coriolis system
 
 As we can observe, there are a few results on the global well-posedness for the tridimensional Boussinesq system $(BQ)^{(\nu_1,\nu_2,\dots,\nu_N)}_{(\eta_1,\eta_2,\dots,\eta_N)}$, when $\nu_i\geq0$ and $\eta_i\geq 0$, for each $1\leq i\leq N$. Most of the global existence results are given considering a non-constant external force multiplying the term $\theta$ at the right hand of the system, or some geometrical structure for the initial data, for instance, axisymmetric initial data. Our global well-posedness result allows to consider initial data in Fourier-Besov spaces, which, for appropiated indexes contain homogeneous functions with negative degree and also large initial data in the $L^{2}$ - space. In order to obtain this result, we use a called rescaled approach for the system in order to deal with the linear force (a constant force multiplying the temperature) at the right hand of the system by mean of the fixed point lemma (Lemma \ref{lema4.1}). Moreover this rescaled approach allows to describe a qualitative behaviour of the system, for instance, we can consider large norm for the initial temperature, if we consider small viscosity and large diffusivity, and we can take large norm for both the initial temperature and velocity, if we consider small diffusivity and large viscosity.
 
%%%%%%%%%%%%%%%%%%%%%%%%%%%%%%
%%%%%%%%%%%%%%%%%%%%%%%%%%%%%%%%%%%%%%%%%%
%%%%%%%%%%%%%%%%%%%

The structure of this paper is the following: In Section \ref{preliminar} we introduce some analysis tools in order to establish the main language for this work, for instance, the Fourier-Besov spaces, the notion of a rescaled Boussinesq system and a fixed point lemma we shall use in the proof of the main theorem. In section \ref{main} it is established the main theorem. In section \ref{estimatesWellposedness} are proved both the estimatives for the linear and the bilinear terms. Finally, in Section \ref{lastsection} is proved the main theorem and it is analised the behaviour of the Boussinesq system with respect to the relation between both the parameters of viscosity and diffusivity and the initial data.

%%%%%%%%%%%%%%%%%%%%%%%%%%%%%%%%%
%%%%%%%%%%%%%%%%%%%%%%%%%%%%%%%%
%%%%%%%%%%%%%%%%%%%%%%%%%%%%%%%%
%%%%%%%%%%%%%%%%%%%%%%%%%%%%%%%%%%%
%%%%%%%%%%%%%%%%%%%%%%%%%%%%%%%%%%%%%%
%%%%%%%%%%%%%%%%%%%%%%%%%%%%%%
\section{Preliminaries}\label{preliminar}

In this section, we recall the definition of Fourier-Besov spaces, we introduce our concept of a rescaled Boussinesq-Coriolis system, we recall the explicit expression for both the semigroups, the heat semigroup and the Stokes-Coriolis semigroup, and the notion of mild solutions. We also recall the paraproduct formula, an abstract fixed point lemma and time-dependent spaces, which are related to our solution spaces.

\subsection{Fourier-Besov spaces}
Fourier-Besov spaces were introduced in the work of P. Konieczny and T. Yoneda (\cite{Konieczny-Yoneda2011}, 2011) to study the dispersive effect of the Coriolis force for both the stationary and the non-stationary Navier-Stokes system. 

In order to introduce these spaces, they consider a radially symmetric function $\varphi\in {\cal S}(\mathbb{R}^{3})$ supported in the annulus $\{\xi\in\mathbb{R}^{3};\frac{3}{4}\leq \mid \xi\mid \leq \frac{8}{3}\}$ such that 

\begin{equation*}
\displaystyle{\sum_{j\in\mathbb{Z}}}\,\varphi(2^{-j}\xi)=1,\,\,\mbox{for all}\,\,\xi\neq 0,    
\end{equation*}
 which permits to introduce two associated functions
 
\begin{equation*}
    \varphi_j(\xi)=\varphi(2^{-j}\xi)\,\,\mbox{and}\,\,\psi_j(\xi)=\displaystyle{\sum_{k\leq j-1}}\varphi_k(\xi)
\end{equation*}
and then to introduce the standard localization operators:

\begin{equation*}
\Delta_jf=\varphi_j(D)f,\,\,\,S_jf=\displaystyle{k\leq j-1}\Delta_kf=\psi_j(D)f,\,\,\,\mbox{for every}\,\,j\in\mathbb{Z}.
\end{equation*}
These operators satisfy the following properties:

\begin{equation*}
\Delta_j\Delta_kf=0,\,\,\mbox{if}\,\mid j-k\mid \geq 2,\,\mbox{and}    
\end{equation*}

\begin{equation*}
\Delta_j(S_{k-1}f \Delta_kf)=0,\,\,\mbox{if}\,\mid j-k\mid \geq 5.    
\end{equation*}
\begin{obs}
The space $\dot{{\cal FB}}_{p,q}^{s}$ contains homogeneous
functions of degree $-d=s-n+\frac{n}{p}$.
\end{obs}
They introduced the called homogeneus Fourier-Besov spaces.

\begin{defi}

    \begin{itemize}
        \item [I.]
       For $1\leq p\leq \infty$, $1\leq q< \infty$ and $s\in \mathbb{R}$,
       
    \begin{equation*}
    \dot{ {\cal FB}}^{s}_{p,q}(\mathbb{R}^{n})=\{
    f\in {\cal S}': \hat{f}\in L^{1}_{loc}, \parallel f\parallel_{\dot{{\cal FB}}^{s}_{p,q}(\mathbb{R}^{n})}=
    ( \displaystyle{\sum_{k\in\mathbb{Z}}}2^{ksq}\parallel \varphi_k\hat{f}\parallel_{L^{p}(\mathbb{R}^{n})}^{q} )^{1/q}<+\infty
    \}
    \end{equation*}
     \item [II.] For $1\leq p\leq \infty$ and $s\in \mathbb{R}$, 
\begin{equation*}
    \dot{{\cal FB}}^{s}_{p,\infty}(\mathbb{R}^{n})=\{
    f\in {\cal S}': \hat{f}\in L^{1}_{loc}, \parallel f\parallel_{\dot{{\cal FB}}^{s}_{p,\infty}(\mathbb{R}^{n})}=\displaystyle{\sup_{k\in\mathbb{Z}}}2^{ks}\parallel\varphi_k\hat{f}\parallel_{L^{p}(\mathbb{R}^{n})}<+\infty\}
\end{equation*}
\end{itemize}
\end{defi}
The next lemma describes some properties related to Fourier-Besov spaces, see Section 2.1 from (\cite{Almeida-Ferreira-Lima-UniformNSCFBMorrey}, 2017) in the case $\nu=0$ and $\nu_1=\nu_2=0$, which corresponds to Fourier-Besov spaces.

\begin{lema}
Let $s_{1},s_{2}\in\mathbb{R}$, $1\leq p_{1},p_{2}<+\infty$.

\begin{itemize}

\item[(i)] (Bernstein-type inequality) Let $p_{2}\leq p_{1}$. If $A>0$ and
$\mbox{supp}\,(\hat{f})\subset\{\xi\in\mathbb{R}^{n};\mid\xi\mid\leq A2^{j}\}$, then

\begin{equation*}
\parallel\xi^{\beta}\hat{f}\parallel_{p_{2}}\leq C2^{j\mid\beta
\mid+j(\frac{n}{p_{2}}-\frac{n}{p_{1}})}\parallel\hat
{f}\parallel_{p_{1}},
\end{equation*}
where $\beta$ is the multi-index, $j\in\mathbb{Z}$, and $C>0$ is a constant
independent of $j,\xi$ and $f$.

\item[(ii)] (Sobolev-type embedding) For $p_{2}\leq p_{1}$ and $s_{2}\leq
s_{1}$ satisfying $s_{2}+\frac{n}{p_{2}}=s_{1}+\frac{n}{p_{1}%
}$, we have the continuous inclusion

\begin{equation*}
\dot{{\cal FB}}_{p_{1},r_{1}}^{s_{1}}\subset\dot{{\cal FB}}_{p_{2},r_{2}}^{s_{2}},
\end{equation*}
for all $1\leq r_{1}\leq r_{2}\leq\infty$.
\end{itemize}
\end{lema}
%%%%%%%%%%%%%%%%%%%%%%%%%%%%%%%%%%%%%%%%%%%%%%%%%%%
%%%%%%%%%%%%%%%%%%%%%%%%%%%%%%%%%%%%%%%%%%%%%%%%%%%

\subsection{The rescaled Boussinesq(-Coriolis) system}
Since the approach we deal with, that is, the search for global mild solutions to the Boussinesq system  $(BQ)^{(\nu,\nu,\nu)}_{(\eta,\eta,\eta)}$, with $\nu>0$ and $\eta>0$, in Fourier-Besov spaces, does not have additional difficulties if we work with the more general Boussinesq-Coriolis system below $(BC)^{(\nu,\nu,\nu)}_{(\eta,\eta,\eta)}$, for $\nu>0$ and $\eta>0$, with Coriolis parameter $\Omega\in \mathbb{R}$, because the obtained results are uniform with respect to the Coriolis parameter $\Omega$; from now on we consider the next system

\begin{equation*}
(BC)^{(\nu,\nu,\nu)}_{(\eta,\eta,\eta)}\,\,\,\left\{
\begin{split}
&  \partial_{t}u-\nu\Delta u+\Omega e_{3}\times u+ (u\cdot\nabla) u +\nabla P=g\theta e_3,\\
&  \partial_{t}\theta-\eta\Delta \theta+(u\cdot\nabla) \theta =0,\\
&  \mbox{div}\,u=0,\,\,\mbox{for}\,\,(x,t)\in\mathbb{R}^{3}\times(0,\infty
),\,\,\mbox{and}\\
&  u(x,0)=u_{0}(x),\,\,\theta(x,0)=\theta_{0}(x),\mbox{for}\,\,x\in
\mathbb{R}^{3},
\end{split}
\right.  %\label{rescaledBQ}
\end{equation*}
which, in the particular case $\Omega=0$ and $g=1$, corresponds to the Boussinesq system $(BQ)^{(\nu,\nu,\nu)}_{(\eta,\eta,\eta)}$, with $\nu>0$ and $\eta>0$.

The main idea for the approach in this work is the notion of rescaled Boussinesq-Coriolis system. For a positive parameter $\lambda$, let us consider the rescaled variables $v=u/\lambda$ and $\rho=\theta/\lambda^{2}$; these new variables satisfy the following rescaled Boussinesq-Coriolis system, 

\begin{equation}
(BC_{\lambda})^{(\nu,\nu,\nu)}_{(\eta,\eta,\eta)}\,\,\,\left\{
\begin{split}
&  \partial_{t}v-\nu\Delta v+\Omega e_{3}\times v+\lambda (v\cdot\nabla) v +\nabla q=\lambda g\rho e_3,\\
&  \partial_{t}\rho-\eta\Delta \rho+\lambda(v\cdot\nabla) \rho =0,\\
&  \mbox{div}\,v=0,\,\,\mbox{for}\,\,(x,t)\in\mathbb{R}^{3}\times(0,\infty
),\,\,\mbox{and}\\
&  v(x,0)=v_{0}(x), \rho(x,0)=\rho_{0}(x),\,\mbox{for}\,\,x\in
\mathbb{R}^{3},
\end{split}
\right.  \label{rescaledBQ}%
\end{equation}
where $q=P/\lambda$ corresponds to the rescaled pressure term.
This kind of rescaling for the original variables appear to be useful in order to obtain the norm of the linear term $\lambda g\rho e_3$ less than one and the bilinear term $\lambda (v\cdot\nabla) v $ small enough, for small $\lambda$, in order to get a global well-posedness result by the application of the fixed point lemma, Lemma \ref{lema4.1}.
%%%%%%%%%%%%%%%%%%%
%%%%%%%%%%%%%%%%%%%%%%%%%%

\subsection{Heat semigroup, Stokes-Coriolis semigroup and mild solutions}
In this section, we recall both the semigroups, the heat semigroup, denoted by $S$, which is given in Fourier variables by the expression

\begin{equation}\label{heatSemigroup}
 (S\rho_0)^{\wedge}(t,\xi)=e^{-\eta\mid\xi\mid^{2}t}\rho_0(\xi)   
\end{equation}
 and the Stokes-Coriolis semigroup, denoted by $S_{\Omega}$, which in Fourier variables is given by
 
\begin{equation}\label{StokesCoriolisSemigroup}
(S_{\Omega}v_0)^{\wedge}(t,\xi)=\cos(\Omega\frac{\xi_3}{\mid\xi\mid})e^{-\nu\mid\xi\mid^{2}t}I\hat{v}_0(\xi)+
\sin(\Omega\frac{\xi_3}{\mid\xi\mid})e^{-\nu\mid\xi\mid^{2}t}R(\xi)\hat{v}_0(\xi),
\end{equation}
where $t\geq 0$, $\xi\in \mathbb{R}^{3}$, $I$ is the identity matrix and the matrix $R(\xi)$ is given by

\begin{equation*}
R(\xi)=\left(
\begin{matrix}
0 & \frac{\xi_3}{\mid\xi\mid} & -\frac{\xi_2}{\mid\xi\mid}\\
-\frac{\xi_3}{\mid\xi\mid} &  0 & \frac{\xi_1}{\mid\xi\mid}\\
\frac{\xi_2}{\mid\xi\mid} & -\frac{\xi_1}{\mid\xi\mid} & 0
\end{matrix}
\right). 
\end{equation*}
It is important to observe that the norm of the matrix $\parallel R(\xi)\parallel\leq 2$, for all $\xi\in\mathbb{R}^{3}$.
Let us redefine the variables $V=v$ and $D=\rho$. Now, if we consider the Duhammel's principle, the system $(BC_{\lambda})^{(\nu,\nu,\nu)}_{(\eta,\eta,\eta)}$ can be rewritten as the integral equation given by 

 \begin{equation}\label{rescaledintegralequation1}
 \left\lbrace
 	\begin{split}
 &V= \frak{V}_0 + L_1(V,D)+B^{1}((V,D), (V,D) )\\
 &D = \frak{D} _0 + L_2(V,D)+B^{2}((V,D), (V,D) ) ,
 \end{split}
 \right.
 \end{equation}  
with 

\begin{equation*}
\begin{split}
&\frak{V}_{0}(t)=S_\Omega(t) v_0\\
&\frak{D}_0(t)=S(t)\rho_0 \\
&L_1(V,D)(t)=  \lambda g\displaystyle{\int_{0}^{t}}S_{\Omega}(t-t')D(t')e_3\,dt' 
\\
&L_2(V,D)(t)= 0 \\
& B^{1}((V,D),(V',D') )(t) = -\frac{1}{2}\lambda \displaystyle{\int_{0}^{t}}S_{\Omega}(t-t')\, \mbox{div} (V(t')\otimes V'(t') )\,dt'-\frac{1}{2}\lambda \displaystyle{\int_{0}^{t}}S_{\Omega}(t-t')\, \mbox{div} (V'(t')\otimes V(t'))\,dt'\\
&\hspace{1.65in}= B^{1,a}((V,D),(V',D') )(t) + B^{1,b}((V,D),
(V',D') )(t)\,\,\,\,\mbox{and}\\
& B^{2}((V,D),(V',D') )(t) =
- \frac{1}{2}\lambda\displaystyle{\int_{0}^{t}}S(t-t') \mbox{div} (V(t')\otimes D'(t') )\,dt'
- \frac{1}{2}\lambda\displaystyle{\int_{0}^{t}}S(t-t') \mbox{div} ( V'(t')\otimes D(t'))\,dt'\\
&\hspace{1.65in}= B^{2,a}((V,D),(V',D') )(t) + B^{2,b}((V,D),
(V',D') )(t).\end{split}
\end{equation*}
A solution in an appropiated Banach space $Z$ for the integral equation $(\ref{rescaledintegralequation1})$ is called a mild solution for the rescaled Boussinesq-Coriolis system $(BC_{\lambda})^{(\nu,\nu,\nu)}_{(\eta,\eta,\eta)}$.
\begin{obs}\label{remarkLB}
If we denote by $X=(V,D)$, $X'=(V',D')$, $X_0=(\frak{V}_{0}(\cdot),\frak{D}_{0}(\cdot))$, $L(X)=(L_1(X),L_2(X))$ and $B(X,X')=(B^{1}(X,X'),B^{2}(X,X'))$ the above system (\ref{rescaledintegralequation1}) can be rewritten as

\begin{equation}\label{rescaledintegralequation2}
    X=X_0+L(X)+B(X,X).
\end{equation}
Therefore, in order to apply Lemma \ref{lema4.1} below, if we denote $Z=\cal{X}\times \cal{Y}$ as in Section \ref{main}, we should to consider that, if
the map $B$ is bilinear from $Z\times Z$ to $Z$ and the operator $L$ is a continuous linear map on $Z$ with norm less than one, we just need to consider small enough initial data $\parallel X_0\parallel_{Z}$ in order to get a unique solution in a small enough ball in $Z$ centered in zero.
\end{obs}

\subsection{Paraproduct formula, abstract f{}ixed point lemma and
time-dependent spaces}
In this subsection we set the Bony's paraproduct formula, describe an abstract
f{}ixed point lemma and we recall suitable time-dependent spaces based on Fourier-Besov spaces which are the required solution spaces. 

For two given tempered distributions $f,g\in\mathcal{S}^{\prime}(\mathbb{R}^{n})$, it is defined  the Bony's
paraproduct operator $T_{f}(\cdot)$ as the following expression

\begin{equation*}
T_{f}(g)=\displaystyle{\sum_{j\in\mathbb{Z}}}S_{j-1}f\Delta_{j}g, \label{Top}%
\end{equation*}
also the operator $R(\cdot,\cdot)$ is defined by 

\begin{equation*}
R(f,g)=\displaystyle{\sum_{j\in\mathbb{Z}}}\Delta_{j}f\tilde{\Delta}%
_{j}g,\,\,\,\mbox{with}\,\,\,\tilde{\Delta}_{j}g=\displaystyle{\sum_{\mid
j^{\prime}-j\mid\leq1}}\Delta_{j^{\prime}}g. \label{Top2}%
\end{equation*}
The symbol of the operator $\tilde{\Delta}_{j}$ in (\ref{Top2}) is denoted by
$\tilde{\varphi}_{j}$. By mean of these operators, we can express the paraproduct of $f$ by $g$ as

\begin{equation}\label{paraproduct}
fg=T_{f}(g)+T_{g}(f)+R(f,g). %
\end{equation}
Expression (\ref{paraproduct}) is known as Bony's paraproduct formula, see (\cite{BonyParaproduct}, 1981) and (\cite{Lemarie2002}, 2002). \newline
%%%%%%%%%%%%%%%%%%%%%
We shall apply the next fixed point lemma in a general Banach space $Z$.
\begin{lema}\label{lema4.1}
	Let be $Z$ a Banach space, $L$ a continuous linear map from $Z$ to $Z$ 
	and let be $B$ a bilinear map from $Z\times Z$ to $Z$.
	Let us def{}ine 
	
	\begin{equation*}
	\vv L\vv_{\mathfrak{L}(Z)}=\displaystyle{\sup_{\vv x\vv=1}}\vv Lx\vv\,\mbox{and}\,
	\vv B\vv_{\mathfrak{B}(Z)}= \displaystyle{\sup_{\vv x\vv=\vv y\vv=1}}
	\vv B(x,y)\vv.
	\end{equation*}
	If $\vv L\vv_{\mathfrak{L}(Z)} < 1$ then for any $x_0\in Z$ such that
	
	\begin{equation*}
	\vv x_0\vv_Z<\frac{(1-\vv L\vv_{\mathfrak{L}(Z)  })^{2}}{4\vv B\vv_{\mathfrak{B}(Z)}},
	\end{equation*}	
	the equation 
	
	\begin{equation*}
	x=x_0+Lx+B(x,x)
	\end{equation*}
	has unique solution in the ball of center $0$ and radius 
	$\frac{1- \vv L\vv_{\mathfrak{L}(Z)}}{2\vv B\vv_{\mathfrak{B}(Z)} }$.
\end{lema}
%%%%%%%%%%%%%%%%%%%%%%%

Finally, we recall two time-dependent spaces. Let be $1\leq p,q\leq\infty$, $1\leq r<\infty$, 
$0<T\leq\infty$ and denote $I=(0,T)$. The Banach spaces $L^{p}(I;\mathcal{FB}%
_{q,r}^{s})$ and $\mathcal{L}^{p}(I;\mathcal{FB}_{q,r}^{s})$ are the
set of Bochner measurable functions from $I$ to $\mathcal{FB}_{q,r}^{s}$
with respective norms given by

\begin{equation*}
\parallel f\parallel_{L^{p}(I;\mathcal{FB}_{q,r}^{s})}%
=\displaystyle{\parallel\parallel f(\cdot,t)\parallel_{\mathcal{FB}_{q,r}^{s}}\parallel_{L^{p}(I)}}=\displaystyle{\parallel
\displaystyle{\left(\displaystyle{\sum_{j\in\mathbb{Z}}}\left(2^{js}\parallel
\varphi_{j}\hat{f}\parallel_{L^{q}(\mathbb{R}^{3})}\right)^{r}\right)^{1/r}\parallel_{L^{p}(I)}}}%
\end{equation*}
and

\begin{equation*}
\parallel f\parallel_{\mathcal{L}^{p}(I;\mathcal{FB}_{q,r}^{s}%
)}=\displaystyle{\left(\displaystyle{\sum_{j\in\mathbb{Z}}\left(2^{js}%
\parallel\varphi_{j}\hat{f}\parallel_{L^{p}(I;L^{q}(\mathbb{R}^{3}))}\right)^{r}%
}\right)^{1/r}}=\left(\displaystyle{\sum_{j\in\mathbb{Z}}}%
\left(2^{js}\parallel\parallel\varphi_{j}\hat{f}\parallel_{L^{q}(\mathbb{R}^{3})}\parallel
_{L^{p}(I)}\right)^{r}\right)^{1/r}.
\end{equation*}
In the case of $r=\infty$, we use the supremum norm in $j\in\mathbb{Z}$. 
%%%%%%%%%%%%%%%%%%
%%%%%%%%%%%%%%%%%%%%%%%%%%%%%%
%%%%%%%%%%%%%%%%%%%%%%%%%%%%%%%%%
\section{Main Result}\label{main}
In this section we establish our main theorem.
 Let us use the notation 
 $\mathcal{X}=
{\cal L}^{\infty}\left(I;\dot{\cal{FB}}^{2-\frac{3}{p}}_{p,q}\right)\cap 
{\cal L}^{1}\left(I;\dot{\cal{FB}}^{4-\frac{3}{p}}_{p,q}\right)$,
${\cal Y}=
{\cal L}^{\infty}\left(I;\dot{\cal{FB}}^{-\frac{3}{p}}_{p,q}\right)\cap 
{\cal L}^{1}\left(I;\dot{\cal{FB}}^{2-\frac{3}{p}}_{p,q}\right)$ and $Z={\cal X}\times{\cal Y} $, where $I$ denotes the interval $(0,\infty)$.
\newpage
\begin{teo}\label{teorema}
Let be $p$ and $q$ two parameters satisfying one of the following conditions:
\begin{itemize}
\item [($c_1$)] $3< p\leq \infty$ and $1\leq q\leq \infty$, or
\item [($c_2$)] $p=3$ and $q=1$.
\end{itemize} 
Then, for all Coriolis parameter $\Omega\in \mathbb{R}$, there exist two positive real numbers $\lambda_0$, depending of $\nu$, and $\varepsilon_0=\varepsilon_0(\lambda_0)$ such that, for values of $\lambda\in\mathbb{R}$ with $0<\lambda\leq \lambda_0$ it is valid the following implication:
If $(V_0,D_0)\in \dot{\cal{FB}}^{2-\frac{3}{p}}_{p,q}\times\dot{\cal{FB}}^{-\frac{3}{p}}_{p,q}$ satisfies

\begin{equation*}
\frac{1}{\nu}\parallel V_0\parallel_{\dot{{\cal FB}}^{2-\frac{3}{p}}_{p,q}}+   
\frac{1}{\eta}\parallel D_0\parallel_{\dot{{\cal FB}}^{-\frac{3}{p}}_{p,q}}<\frac{\varepsilon_0}{C_0},
\end{equation*}
for some constant $C_0>0$, then there exist a mild solution $(V,D)\in Z$ for the Boussinesq-Coriolis system $(BC_{\lambda})^{(\nu,\nu,\nu)}_{(\eta,\eta,\eta)}$ and this solutions is unique in the ball of center $0$ and radius $\frac{1-\parallel L\parallel_{{\cal L}(Z)}}{2\parallel B\parallel_{{\cal B}(Z)}}$.    
\end{teo}

\begin{obs}
It is important to note that the parameter $\lambda_0$ does not depends on the Coriolis parameter $\Omega$, so that this result is also valid for the rescaled Boussinesq system $(BQ_{\lambda})^{(\nu,\nu,\nu)}_{(\eta,\eta,\eta)}$, with $\nu>0$ and $\eta>0$, which corresponds to  $(BC_{\lambda})^{(\nu,\nu,\nu)}_{(\eta,\eta,\eta)}$ with $\Omega=0$.

\end{obs}

%%%%%%%%%%%%%%%%%%%%%%%%%%%%%%
%%%%%%%%%%%%%%%%%%%%%%%%%%%
%%%%%%%%%%%%%%%%%%%%%%%%%%%
%%%%%%%%%%%%%%%%%%%%%%
\section{Estimates for the well-posedness of the rescaled Boussinesq-Coriolis system}\label{estimatesWellposedness}
In this section we shall prove estimates for the linear operator $L$, the bilinear map $B$, as in Remark \ref{remarkLB}, in order to apply directly Lemma \ref{lema4.1}.
Let us start with estimates for the linear term $L=(L_1,L_2)$.
\begin{lema}
Let $1\leq p\leq \infty$ and $1\leq q\leq \infty$ be two parameters. For all Coriolis parameter $\Omega\in\mathbb{R}$, the following estimate it is valid for the linear operator $L$:

\begin{equation*}
\parallel L\parallel_{{\cal L}(Z)}\leq \frac{C_L\lambda}{\nu}\mid g\mid,    
\end{equation*}
where $C_{L}$>0 is a constant independent on $\Omega$.
\end{lema}
 Since $L_2\equiv 0$, we just need to estimate the linear operator $L_1$. Recal that $X=(V,D)$, see Remark \ref{remarkLB}. By definition of the ${\cal X}$-norm, we estimate
 
\begin{equation}\label{A1}
    \parallel L_1X\parallel_{{\cal X}}\leq \lambda \mid g\mid ( I_1+I_2),
\end{equation}
where 

\begin{equation*}
I_1=\displaystyle{\parallel}
2^{j(2-\frac{3}{p})} \displaystyle\parallel
\displaystyle{\int_{0}^{t}}e^{-\nu(t-t')2^{2j}}\varphi_j(\xi)\hat{D}(\xi)\displaystyle{\parallel}_{L^{\infty}(I;L^{p})}
\displaystyle{\parallel}_{l^{q}(\mathbb{Z})}
\end{equation*}
and 

\begin{equation*}
I_2=\displaystyle{\parallel}
2^{j(4-\frac{3}{p})} \displaystyle\parallel
\displaystyle{\int_{0}^{t}}e^{-\nu(t-t')2^{2j}}\varphi_j(\xi)\hat{D}(\xi)\displaystyle{\parallel}_{L^{1}(I;L^{p})}
\displaystyle{\parallel}_{l^{q}(\mathbb{Z})}.
\end{equation*}
For the first term $I_1$, by using the Young's inequality and solving the integral

\begin{equation*}
    \displaystyle{\int_{0}^{\infty}}e^{-\nu s2^{2j}}\,ds=\frac{2^{-2j}}{\nu}
\end{equation*}
we have

\begin{equation*}
\begin{split}
I_1\leq &  \parallel 2^{j(2-\frac{3}{p})}\parallel
\displaystyle{\int_{0}^{t}}e^{-\nu(t-t')2^{2j}}\parallel \varphi_j\hat{D}\parallel_{L^{p}}\,dt'\parallel_{L^{\infty}(I)}
\parallel_{l^{q}(\mathbb{Z})}\\
&
\leq C_1
\parallel 2^{j(2-\frac{3}{p})}\parallel e^{-\nu s2^{2j}}\parallel_{L^{1}(I)}\parallel\varphi_j\hat{D}\parallel_{L^{\infty}(I;L^{p})}
\parallel_{l^{q}(\mathbb{Z})}\\
&
\leq \frac{C_1}{\nu}
\parallel D\parallel_{{\cal L}^{\infty}(I;
\dot{{\cal FB}}^{-\frac{3}{p}}_{p,q})}
\end{split}
\end{equation*}
and thus we get

\begin{equation}\label{A2}
I_1\leq\frac{C_1}{\nu}\parallel X\parallel_Z.    
\end{equation}
Similarly, by using the Young's inequality, we have

\begin{equation*}
\begin{split}
I_2 &\leq \parallel 2^{j(4-\frac{3}{p})} \parallel 
\displaystyle{\int_{0}^{t}}e^{-\nu(t-t')2^{2j}}
\parallel \varphi_j\hat{D}\parallel_{L^{p}}\,dt'\parallel_{L^{1}(I)}\parallel_{l^{q}(\mathbb{Z})}\\
&\leq
C_2 \parallel 2^{j(4-\frac{3}{p})}\parallel e^{-\nu s 2^{2j}}\parallel_{L^{1}(I)}\parallel\varphi_j\hat{D}\parallel_{L^{1}(I;L^{p})}\parallel_{l^{q}(\mathbb{Z})}\\
&\leq\frac{C_2}{\nu}\parallel D\parallel_{{\cal Y}}
\end{split}
\end{equation*}
and thus

\begin{equation}\label{A3}
I_2\leq \frac{C_2}{\nu}\parallel X\parallel_{Z}.    
\end{equation}
Substituing (\ref{A2}) and (\ref{A3}) in (\ref{A1}), we get

\begin{equation*}
\parallel L_1 X\parallel_{{\cal X}}\leq \frac{C_{L}\lambda}{\nu} \mid g\mid \parallel X\parallel_Z     
\end{equation*}
and therefore, since $L_2\equiv 0$, we conclude that

\begin{equation}\label{A4}
    \parallel L\parallel_{{\cal L}(Z)}\leq \frac{C_L\lambda}{\nu}\mid g\mid.
\end{equation}
%%%%%%%%%%%%%%%%%%%%%%%%%%%%%%%%%%%%%%%%%
For the estimate of the bilinear map $B_2$, we establish the next lemma. 
\begin{lema}\label{lemaSemigroup}
For $1\leq p\leq \infty$ and $1\leq q\leq \infty$, are valid the following items (i) and (ii) below:
\begin{itemize}
    \item [(i)]
    
\begin{equation}\label{A6}
\parallel \displaystyle{\int_{0}^{t}}S(t-t')f(t')\,dt'
\parallel_{{\cal L}^{\infty}(I;\dot{\cal{FB}}^{-\frac{3}{p}}_{p,q})}\leq C
\parallel f\parallel_{{\cal L}^{1}(I;\dot{\cal{FB}}^{-\frac{3}{p}}_{p,q})}
\end{equation}
\item [(ii)]

\begin{equation}\label{A7}
\parallel \displaystyle{\int_{0}^{t}}S(t-t')f(t')\,dt'
\parallel_{{\cal L}^{1}(I;\dot{\cal{FB}}^{2-\frac{3}{p}}_{p,q})}\leq \frac{C}{\eta}
\parallel f\parallel_{{\cal L}^{1}(I;\dot{\cal{FB}}^{-\frac{3}{p}}_{p,q})},
\end{equation}
\end{itemize}
for some positive constant $C$.
\end{lema}
Let us start with the estimate of item (i): From definition of the space ${\cal L}^{\infty}(I;\dot{\cal{FB}}^{-\frac{3}{p}}_{p,q})$ and applying the Young's inequality, we have the estimates

\begin{equation*}
\begin{split}
\parallel \displaystyle{\int_{0}^{t}}S(t-t')f(t')\,dt'
\parallel_{{\cal L}^{\infty}(I;\dot{\cal{FB}}^{-\frac{3}{p}}_{p,q})}&=\parallel 2^{j(-\frac{3}{p})} \parallel 
\displaystyle{\int_{0}^{t}}e^{-\eta(t-t')2^{2j}}\varphi_{j}(\xi)\hat{f}(t',\xi)\,dt'\parallel_{L^{\infty}(I;L^{p})}
\parallel_{l^{q}(\mathbb{Z})}\\
&\leq \parallel 2^{j(2-\frac{3}{p})}\parallel 
\displaystyle{\int_{0}^{t}}e^{-\eta (t-t')2^{2j}}\parallel 
\varphi_j(\xi)\hat{f}(t',\xi)\parallel_{L^{p}}\parallel_{L^{\infty}(I)}
\parallel_{l^{q}(\mathbb{Z})}\\
&\leq C\parallel 2^{j(2-\frac{3}{p})}\parallel e^{-\eta s2^{2j}}\parallel_{L^{\infty}(I)} \parallel \varphi_j(\xi)\hat{f}(t',\xi)\parallel_{L^{1}(I;L^{p})}\parallel_{l^{\mathbb{Z}}}
\end{split}    
\end{equation*}
and, since $\parallel e^{-\eta s2^{2j}}\parallel_{L^{\infty}(I)}\leq 1$, we get the desired estimate (\ref{A6}). Now, we shall prove estimate (\ref{A7}): Again, by definition of the space ${\cal L}^{1}(I;\dot{{\cal FB}}^{2-\frac{3}{p}}_{p,q})$ and Young's inequality, we have the estimates

\begin{equation*}
\begin{split}
\parallel \displaystyle{\int_{0}^{t}}S(t-t')f(t')\,dt'
\parallel_{{\cal L}^{1}(I;\dot{\cal{FB}}^{2-\frac{3}{p}}_{p,q})}
&\leq \parallel 2^{j(2-\frac{3}{p})}\parallel 
\displaystyle{\int_{0}^{t}}e^{-\eta (t-t')2^{2j}}\parallel 
\varphi_j(\xi)\hat{f}(t',\xi)\parallel_{L^{p}}\parallel_{L^{1}(I)}
\parallel_{l^{q}(\mathbb{Z})}\\
&\leq C\parallel 2^{j(2-\frac{3}{p})}\parallel e^{-\eta s2^{2j}}\parallel_{L^{1}(I)} \parallel \varphi_j(\xi)\hat{f}(t',\xi)\parallel_{L^{1}(I;L^{p})}\parallel_{l^{q}(\mathbb{Z})}
\end{split}    
\end{equation*}
and since $\parallel e^{-\eta s2^{2j}}\parallel_{L^{1}(I)}\leq \frac{2^{-2j}}{\eta}$, we get the desired estimate (\ref{A7}).
%%%%%%%%%%%%%%%%%%%%%%%%%%%%%%%%%%%%%
\newline
In the next lemma are summarized bilinear estimates for operators $B_1$ and $B_2$.
\begin{lema}
Let be $p$ and $q$ two parameters satisfying one of the following conditions:
\begin{itemize}
\item [($c_1$)] $3< p\leq \infty$ and $1\leq q\leq \infty$, or
\item [($c_2$)] $p=3$ and $q=1$.
\end{itemize}
Then, for all Coriolis parameter $\Omega\in\mathbb{R}$, the following estimates for the bilinear terms are valid.\\
(i)

\begin{equation}\label{A8}
\parallel B_1(X,X')    
\parallel_{{\cal X}}\leq C_{B_1}\cdot \lambda\cdot\max(1,\frac{1}{\nu}) \parallel X\parallel_{Z}\parallel X'\parallel_{Z}\,\mbox{and}
\end{equation}
(ii)

\begin{equation}\label{A9}
\parallel B_2(X,X')    
\parallel_{{\cal Y}}\leq C_{B_2}\cdot\lambda\cdot \max(1,\frac{1}\eta) \parallel X\parallel_{Z}\parallel X'\parallel_{Z},
\end{equation}
for some positive constants $C_{B_1}$ and $C_{B_2}$ independent of the Coriolis parameter.
\end{lema}
In order to deal with the estimate (\ref{A8}), we follow the same reasoning and computations made in Lemma 4.2 and Lemma 4.3 from (\cite{Aurazo-Ferreira-Global-BC-Stratification-FBM-2021}, 2021).
Now we shall prove estimate (\ref{A9}): If we consider the term $f=\mbox{div}(V\otimes D')$ in Lemma \ref{lemaSemigroup}, we get the estimate

\begin{equation*}
\parallel B_2(X,X')\parallel_{{\cal Y}}\leq C\lambda \max (1,\frac{1}\eta)\parallel \mbox{div} (V\otimes D')\parallel_{
{\cal L}^{1}(I;\dot{{\cal FB}}^{-\frac{3}{p}}_{p,q})}  
\end{equation*}
and thus we conclude that

\begin{equation}\label{A9'}
\parallel B_2(X,X')\parallel_{{\cal Y}}\leq C\lambda \max (1,\frac{1}\eta)\parallel V\otimes D'\parallel_{
{\cal L}^{1}(I;\dot{{\cal FB}}^{1-\frac{3}{p}}_{p,q})}.
\end{equation}
Therefore, if we show the following estimate

\begin{equation}\label{A10}
\parallel V\otimes D'\parallel_{
{\cal L}^{1}(I;\dot{{\cal FB}}^{1-\frac{3}{p}}_{p,q})}\leq C
\parallel V\parallel_{{\cal X}}\parallel D'\parallel_{{\cal Y}},
\end{equation}
we can conclude the desired estimate (\ref{A9}).
In the rest of the proof we show estimative (\ref{A10}).
By the paraproduct formula (2.4) from (\cite{Konieczny-Yoneda2011}, 2011), we have

\begin{equation}\label{A11}
\varphi_j(VD')^{\wedge}=I^{j}_{1}+I^{j}_{2}+I^{j}_{3},    
\end{equation}
where 

\begin{equation*}
I^{j}_{1}=\displaystyle{\sum_{\mid k-j\mid\leq 4}}\varphi_j
[(S_{k-1}V)^{\wedge}\ast(\varphi_{k}\hat{D'})],
\end{equation*}

\begin{equation*}
I^{j}_{2}=\displaystyle{\sum_{\mid k-j\mid\leq 4}}\varphi_j
[(S_{k-1}D')^{\wedge}\ast(\varphi_{k}\hat{V})]
\end{equation*}
and 

\begin{equation*}
I^{j}_{3}=
\displaystyle{\sum_{k\geq j-2}}\varphi_{j}
[(\varphi_{k}\hat{V})\ast(\tilde{\varphi}_{k}\hat{D'})],
\end{equation*}
where 

\begin{equation*}
\tilde{\varphi}_{k}= \varphi_{k-1}+\varphi_{k}+\varphi_{k+1}.    
\end{equation*}
Let us start with estimates for the term $I^{j}_{1}$: By applying the Young's inequality in $\mathbb{R}^{3}$ we have

\begin{equation*}
\begin{split}
\parallel I^{j}_{1}\parallel_{L^{1}(I;L^{p})}&\leq 
C\displaystyle{\sum_{\mid k-j\mid\leq 4}}
\parallel \parallel (S_{k-1}V)^{\wedge}\parallel_{L^{1}}
\parallel \varphi_{k}\hat{D'}\parallel_{L^{p}}
\parallel_{L^{1}(I)}\\
&\leq C
\displaystyle{\sum_{\mid k-j\mid\leq 4}}\,
\displaystyle{\sum_{k'<k}}2^{k'(3-\frac{3}{p})}
\parallel\varphi_{k'}\hat{V}\parallel_{L^{\infty}(I;L^{p})}\parallel\varphi_{k}\hat{D'}\parallel_{L^{1}(I;L^{p})}
\end{split}    
\end{equation*}
and applying the H{\"o}lder's inequality for series we get, for $1\leq q\leq \infty$, that

\begin{equation*}
\begin{split}
\parallel I^{j}_{1}\parallel_{L^{1}(I;L^{p})}&\leq 
C
\displaystyle{\sum_{\mid k-j\mid\leq 4}}\,
\displaystyle{\sum_{k'<k}}2^{k'} 2^{k'(2-\frac{3}{p})}
\parallel\varphi_{k'}\hat{V}\parallel_{L^{\infty}(I;L^{p})}\parallel\varphi_{k}\hat{D'}\parallel_{L^{1}(I;L^{p})}\\
&\leq C
\displaystyle{\sum_{\mid k-j\mid\leq 4}}\,2^{k}
\parallel\varphi_{k}\hat{D'}\parallel_{L^{1}(I;L^{p})}
\parallel V\parallel_{{\cal L}^{\infty}(I;{\dot{{\cal FB}}^{2-\frac{3}{p}}_{p,q}})}
\end{split}
\end{equation*}
and thus

\begin{equation*}
\begin{split}
2^{j(1-\frac{3}{p})}\parallel I^{j}_{1}\parallel_{L^{1}(I;L^{p})}&\leq 
C\displaystyle_{\sum_{k}}2^{(j-k)(1-\frac{3}{p})}\chi_{\{l;\mid l\mid\leq 4\}}(j-k)2^{k(2-\frac{3}{p})}
\parallel\varphi_{k}\hat{D}\parallel_{L^{1}(I;L^{p})}
\parallel V\parallel_{{\cal L}^{\infty}(I;\dot{{\cal FB}}^{2-\frac{3}{p}}_{p,q})}\\
&\leq C (a_l\ast b_k)_j\parallel V\parallel_{{\cal L}^{\infty}(I;\dot{{\cal FB}}^{2-\frac{3}{p}}_{p,q})},
\end{split}
\end{equation*}
where $a_l=2^{l(1-\frac{3}{p})}\chi_{\{l;\mid l\mid\leq 4\}}(l)$ and $b_k=2^{k(2-\frac{3}{p})}\parallel \varphi_j\hat{D}\parallel_{L^{1}(I;L^{p})}$.
Finally, applying the $l^{q}(\mathbb{Z})$-norm and Young's inequality for series we get

\begin{equation}\label{A14}
\parallel 2^{j(1-\frac{3}{p})}\parallel I^{j}_1\parallel_{L^{1}(I;L^{p})}\parallel_{l^{q}(\mathbb{Z})}
\leq C\parallel D'\parallel_{{\cal L}^{1}(I;\dot{{\cal FB}}^{2-\frac{3}{p}}_{p,q})}\parallel V\parallel_{{\cal L}^{\infty}(I;\dot{{\cal FB}}^{2-\frac{3}{p}}_{p,q})}
\end{equation}
Now, we estimate the term $I^{j}_{2}$: By applying the Young's inequality in $\mathbb{R}^{3}$ and H{\"o}lder's inequality for series we have

\begin{equation*}
\begin{split}
\parallel I^{j}_2\parallel_{L^{1}(I;L^{p})} &\leq 
C\displaystyle{\sum_{\mid k-j\mid\leq 4}}\parallel \parallel 
(S_{k-1}D')^{\wedge}\parallel_{L^{1}} \parallel\varphi_k\hat{V}\parallel_{L^{p}}\parallel_{L^{1}(I)}\\
&\leq
C \displaystyle{\sum_{\mid k-j\mid\leq 4}}\,
\displaystyle{\sum_{k'<k}}2^{k'(3-\frac{3}{p})}
\parallel\varphi_{k'}\hat{D'}\parallel_{L^{\infty}(I;L^{p})}
\parallel\varphi_k\hat{V}\parallel_{L^{1}(I;L^{p})}\\
&\leq 
C \displaystyle{\sum_{\mid k-j\mid\leq 4}}\,
\displaystyle{\sum_{k'<k}}2^{3k'}\cdot2^{k'(-\frac{3}{p})}
\parallel\varphi_{k'}\hat{D'}\parallel_{L^{\infty}(I;L^{p})}
\parallel\varphi_k\hat{V}\parallel_{L^{1}(I;L^{p})}\\
&\leq C
\displaystyle{\sum_{\mid k-j\mid\leq 4}}2^{3k}
\parallel\varphi_k \hat{V}\parallel_{L^{1}(I;L^{p})}\parallel D'\parallel_{{\cal L}^{\infty}(I;\dot{{\cal FB}}^{-\frac{3}{p}}_{p,q})}\\
&\leq C
\displaystyle{\sum_{\mid k-j\mid\leq 4}}2^{k(-1+\frac{3}{p})}\cdot 2^{k(4-\frac{3}{p})}
\parallel\varphi_k \hat{V}\parallel_{L^{1}(I;L^{p})}\parallel D'\parallel_{{\cal L}^{\infty}(I;\dot{{\cal FB}}^{-\frac{3}{p}}_{p,q})}
\end{split}
\end{equation*}
and thus

\begin{equation*}
\begin{split}
2^{j(1-\frac{3}{p})}\parallel I^{j}_2\parallel_{L^{1}(I;L^{p})} &\leq C\displaystyle{\sum_{k}}2^{(j-k)(1-\frac{3}{p})}\chi_{\{l;\mid l\mid\leq 4\}}(j-k)2^{k(4-\frac{3}{p})}\parallel \varphi_k\hat{V}\parallel_{L^{1}(I;L^{p})}\parallel D'\parallel_{{\cal L}^{\infty}(I;\dot{{\cal FB}}^{-\frac{3}{p}}_{p,q})}   \\
&\leq C(a_l\ast b_k)_j\parallel D'\parallel_{{\cal L}^{\infty}(I;\dot{{\cal FB}}^{-\frac{3}{p}}_{p,q})}, 
\end{split}
\end{equation*}
where $a_l=2^{l(1-\frac{3}{p})}\chi_{\{l;\mid l\mid\leq 4\}}(l)$ and $b_k=2^{k(4-\frac{3}{p})}\parallel \varphi_k\hat{V}\parallel_{L^{1}(I;L^{p})}$. By applying the Young's inequality for series we obtain

\begin{equation}\label{A17}
\parallel 2^{j(1-\frac{3}{p})}\parallel I^{j}_2 \parallel_{L^{1}(I;L^{p})}\parallel_{l^{q}(\mathbb{Z})}\leq 
C\parallel V\parallel_{{\cal L}^{1}(I;\dot{{\cal FB}}^{4-\frac{3}{p}}_{p,q})}\parallel D'\parallel_{{\cal L}^{\infty}(I;\dot{{\cal FB}}^{-\frac{3}{p}}_{p,q})}.
\end{equation}
Finally we estimate the term $I^{j}_3$: By applying the Young's inequality in $\mathbb{R}^{3}$ and H{\"o}lder's inequality for series, we get

\begin{equation*}
\begin{split}
\parallel I^{j}_3\parallel_{L^{1}(I;L^{p})}&\leq 
C\displaystyle{\sum_{k\geq j-2}} \parallel \varphi_k\hat{V}\parallel_{L^{1}(I;L^{p})}
\displaystyle{\sum_{\mid k-k'\mid\leq 1}}2^{k'(3-\frac{3}{p})}
\parallel\varphi_{k'}\hat{D'}\parallel_{L^{\infty}(I;L^{p})}\\
&\leq C \displaystyle{\sum_{k\geq j-2}} \parallel \varphi_k\hat{V}\parallel_{L^{1}(I;L^{p})}
\displaystyle{\sum_{\mid k-k'\mid\leq 1}}2^{k'3}\cdot2^{-k'(\frac{3}{p})}\parallel\varphi_{k'}\hat{D'}\parallel_{L^{\infty}(I;L^{p})}\\
&\leq C \displaystyle{\sum_{k\geq j-2}} \parallel \varphi_k\hat{V}\parallel_{L^{1}(I;L^{p})}2^{3k}\parallel D\parallel_{{\cal L}^{\infty}(I;\dot{{\cal FB}}^{-\frac{3}{p}}_{p,q})}
\end{split}
\end{equation*}
and thus

\begin{equation*}
\begin{split}
2^{j(1-\frac{3}{p})}\parallel I^{j}_3\parallel_{L^{1}(I;L^{p})}&\leq C\displaystyle{\sum_{k}}\chi_{\{l;l\leq 2\}}(j-k)2^{(j-k)(1-\frac{3}{p})}2^{k(4-\frac{3}{p})}\parallel\varphi_k\hat{V}\parallel_{L^{1}(I;L^{p})}\parallel D\parallel_{{\cal L}^{\infty}(I;\dot{{\cal FB}}^{-\frac{3}{p}}_{p,q})}\\
&\leq C (a_l\ast b_k)_j\parallel D\parallel_{{\cal L}^{\infty}(I;\dot{{\cal FB}}^{-\frac{3}{p}}_{p,q})}, 
\end{split}
\end{equation*}
where $a_l=2^{(1-\frac{3}{p})l}\chi_{\{l;l\leq 2\}}(l)$ and $b_k=2^{k(4-\frac{3}{p})}\parallel\varphi_k\hat{V}\parallel_{L^{1}(I;L^{p})}$. Applying the Young's inequality for series and by using the assumption $p>3$, we have

\begin{equation}\label{A18}
\parallel 2^{j(1-\frac{3}{p})}\parallel I^{j}_3\parallel_{L^{1}(I;L^{p})}\parallel_{l^{q}(\mathbb{Z})}\leq C\parallel V\parallel_{{\cal L}^{1}(I;\dot{{\cal FB}}^{4-\frac{3}{p}}_{p,q})}\parallel D'\parallel_{{\cal L}^{\infty}(I;\dot{{\cal FB}}^{-\frac{3}{p}}_{p,q})}.
\end{equation}
Observe that, if we consider $p=3$ we just need to take $q=1$ and it will be valid the last estimate for $I^{j}_3$.
Considering estimates (\ref{A14}), (\ref{A17}) and (\ref{A18}) in the decomposition (\ref{A11}), we conclude that

\begin{equation*}
\parallel V\otimes D'\parallel_{{\cal L}^{1}(I;\dot{{\cal FB}}^{1-\frac{3}{p}}_{p,q})} \leq C\left(
\parallel V\parallel_{{\cal L}^{\infty}(I;\dot{{\cal FB}}^{2-\frac{3}{p}}_{p,q})} \parallel D'\parallel_{{\cal L}^{1}(I;\dot{{\cal FB}}^{2-\frac{3}{p}}_{p,q})}+ 
\parallel V\parallel_{{\cal L}^{1}(I;\dot{{\cal FB}}^{4-\frac{3}{p}}_{p,q})} \parallel D'\parallel_{{\cal L}^{\infty}(I;\dot{{\cal FB}}^{-\frac{3}{p}}_{p,q})}
\right)     
\end{equation*}
and thus 

\begin{equation}\label{A19}
\parallel V\otimes D'\parallel_{{\cal L}^{1}(I;\dot{{\cal FB}}^{1-\frac{3}{p}}_{p,q})} \leq C \parallel V\parallel_{{\cal X}}\parallel D'\parallel_{{\cal Y}}    
\end{equation}
as required. Estimatives (\ref{A19}) and (\ref{A9'}) implies

\begin{equation}\label{A20}
\parallel B_2(X.X') \parallel_{{\cal Y}}\leq C\lambda\max (1,\frac{1}{\eta}) \parallel V\parallel_{{\cal X}}\parallel D'\parallel_{{\cal Y}}   
\end{equation}
and therefore the desired estimate(\ref{A9}).
%%%%%%%%%%%%%%%%%%%%%%%%%%%%%%%%%%%
%%%%%%%%%%%%%%%%%%%%%%%%%%%%%%%%%
%%%%%%%%%%%%%%%%%%%%%
\begin{section}{Proof of the main theorem and behaviour of the Boussinesq-Coriolis system}\label{lastsection}

\begin{proof}
By the last lemma, we conclude that 

\begin{equation}\label{A21}
\parallel B(X,X')\parallel_Z\leq  C_B\lambda\max (1,\frac{1}{\nu},\frac{1}{\eta}) \parallel X\parallel_{Z}\parallel X'\parallel_{Z},
\end{equation}
for some positive constant $C_B$ independent of the Coriolis parameter $\Omega$. Also, for all Coriolis parameter $\Omega$, 
we estimate

\begin{equation}\label{A22}
\parallel X_0\parallel_Z=\parallel S_\Omega(t)V_0\parallel_{{\cal X}}+\parallel S(t)D_0\parallel_{{\cal Y}}
\leq C_0\left(\frac{1}{\nu}\parallel V_0\parallel_{\dot{{\cal FB}}^{2-\frac{3}{p}}_{p,q}}+
\frac{1}{\eta}\parallel D_0\parallel_{\dot{{\cal FB}}^{-\frac{3}{p}}_{p,q}}\right),
\end{equation}
for some positive constant $C_0$ independent of the Coriolis parameter $\Omega$. With estimates (\ref{A4}), (\ref{A21}) and (\ref{A22}) in mind, we follow the next reasoning: Let us take $\lambda_0=\frac{\nu}{2C_L\mid g\mid }>0$, so that for $0<\lambda\leq \lambda_0$ we have $\parallel L\parallel_{{\cal L}(Z)}<\frac{1}{2}$. By the direct application of the Lemma \ref{lema4.1}, we just need to choose 

\begin{equation}\label{initialOCondition}
\frac{1}{\nu}\parallel V_0\parallel_{\dot{{\cal FB}}^{2-\frac{3}{p}}_{p,q}}+
\frac{1}{\eta}\parallel D_0\parallel_{\dot{{\cal FB}}^{-\frac{3}{p}}_{p,q}}\leq \frac{\varepsilon_0}{C_0},
\end{equation}
where $C_0$ is the constant from (\ref{A22}) and $\varepsilon_0$ is given by the expression 

\begin{equation*}
\varepsilon_{0}=\varepsilon_{0}(\lambda_0)= \frac{(1-\frac{C_L\lambda_{0}}{\nu}\mid g\mid)^{2}}{4C_B\lambda_0\max(1,\frac{1}{\nu},\frac{1}{\eta})}= C\cdot \frac{\min(1,\nu,\eta)}{\nu}>0,
\end{equation*}
where $C=\frac{C_L\mid g\mid}{16C_B}$, in order to prove the existence of a unique mild solution $(V,D)\in Z$ for the Boussinesq-Coriolis system $(BC_{\lambda})^{(\nu,\nu,\nu)}_{(\eta,\eta,\eta)}$, with $\nu>0$ and $\eta>0$, in the ball of center $0$ and radius $\frac{1-\parallel L\parallel_{{\cal L}(Z)}}{2\parallel B\parallel_{{\cal B}(Z)}}$. Therefore we conclude the proof of the main theorem.
\end{proof}

Let us consider the Boussinesq-Coriolis system $(BC)^{(\nu,\nu,\nu)}_{(\eta,\eta,\eta)}$, with $\nu>0$ and $\eta$, and its mild formulation

\begin{equation*}
\begin{split}
&u=\tilde{u}_0+\tilde{L}_1(u,\theta)+\tilde{B}^{1}((u,\theta),(u,\theta))\\
&\theta=\tilde{\theta}_0+\tilde{L}_2(u,\theta)+\tilde{B}^{2}((u,\theta),(u,\theta)),
\end{split}    
\end{equation*}
where 

\begin{equation*}
\begin{split}
&\tilde{u}_0=S_\Omega(t) u_0\\
&\tilde{\theta}_0(t)=S(t)\theta_0 \\
&L_1(u,\theta)(t)=   g\displaystyle{\int_{0}^{t}}S_{\Omega}(t-t')\theta(t')e_3\,dt' 
\\
&\tilde{L}_2(V,D)(t)= 0 \\
& \tilde{B}^{1}((u,\theta),(u',\theta') )(t) = -\frac{1}{2} \displaystyle{\int_{0}^{t}}S_{\Omega}(t-t')\, \mbox{div} (u(t')\otimes u'(t') )\,dt'-\frac{1}{2} \displaystyle{\int_{0}^{t}}S_{\Omega}(t-t')\, \mbox{div} (u'(t')\otimes u(t'))\,dt'\,\,\,\,\mbox{and}\\
& \tilde{B}^{2}((u,\theta),(u',\theta') )(t) =
- \frac{1}{2}\displaystyle{\int_{0}^{t}}S(t-t') \mbox{div} (u(t')\otimes \theta'(t') )\,dt'
- \frac{1}{2}\displaystyle{\int_{0}^{t}}S(t-t') \mbox{div} ( u'(t')\otimes \theta(t'))\,dt'.\\
\end{split}
\end{equation*}
Let be $\lambda_0$, $\varepsilon_0$ and $C_0$ as given in the main theorem (\ref{main}).
For $0<\lambda\leq \lambda_0$, let us take $u_0\in \dot{\cal FB}^{2-\frac{3}{p}}_{p,q}$ and $\theta_0\in \dot{\cal FB}^{-\frac{3}{p}}_{p,q}$ such that

\begin{equation}\label{originalInitialData}
\frac{1}{\nu\lambda}\parallel u_0\parallel_{\dot{\cal FB}^{2-\frac{3}{p}}_{p,q}}+
\frac{1}{\eta\lambda^{2}}\parallel \theta_0\parallel_{\dot{\cal FB}^{-\frac{3}{p}}_{p,q}}<\frac{\varepsilon_0}{C_0},
\end{equation}
there exist a global mild solution $(V,D)\in Z$ for the Boussinesq-Coriolis rescaled system $(BC_{\lambda})^{(\nu,\nu,\nu)}_{(\eta,\eta,\eta)}$, with initial data $V_0=\frac{u_0}{\lambda}$ and $D_0=\frac{\theta_0}{\lambda^{2}}$. Now, if we consider the variables $u=\lambda V$ and $\theta=\lambda^{2}D$ then the pair $(u, \theta)\in Z$ is a global mild solution for the Boussinesq-Coriolis system $(BC)^{(\nu,\nu,\nu)}_{(\eta,\eta,\eta)}$ with given initial data $u_0,\theta_0$.
Moreover, since the solution $(V, D)\in Z$ of the rescaled Boussinesq-Coriolis system is unique in a determined ball in $Z$, for instance with radius $r=\frac{1-\parallel L\parallel_{{\cal L}(Z)}}{2\parallel B\parallel_{{\cal B}(Z)}}$, the ball $B_{Z}(0,r)$ , then the solution $(u,\theta)\in Z$ is unique in the set $\{(u,\theta)\in Z; \frac{1}{\lambda}\parallel u\parallel_{\cal X}+\frac{1}{\lambda^{2}}\parallel\theta\parallel_{{\cal Y}}<r\}$.
 
Now on, we study the behaviour of the Boussinesq-Coriolis system $(BC)^{(\nu,\nu,\nu)}_{(\eta,\eta,\eta)}$ by mean of the analysis of particular cases from the initial conditions (\ref{initialOCondition}), considering the original variables $(u,\theta)$. In this sense, if we consider the initial data $(u_0, \theta_0)$ satisfying the condition (\ref{originalInitialData}), for $\lambda=\lambda_0$, then it is enough to consider initial velocity 

\begin{equation*}
\parallel u_0\parallel_{\dot{{\cal FB}}^{2-\frac{3}{p}}_{p,q}}< \frac{C}{4C_0C_L\mid g\mid}\cdot\nu\cdot\min(1,\nu,\eta) 
\end{equation*}
and the initial temperature disturbance

\begin{equation*}
\parallel \theta_0\parallel_{\dot{{\cal FB}}^{-\frac{3}{p}}_{p,q}}< \frac{C}{8C_0C_{L}^{2}\mid g\mid^{2}}\cdot\nu\cdot\eta\cdot\min(1,\nu,\eta).
\end{equation*}
Let us denote $C_1=\frac{C}{4C_0C_L\mid g\mid}$ and $C_2=\frac{C}{8C_0C_{L}^{2}\mid g\mid^{2}}$. With this notation, we shall consider all possible inequality relations for the terms $\nu$ and $\eta$, in order to consider some values for $\min(1,\nu,\eta)$.
\begin{itemize}
\item [($a_1$)] ($\nu=1$ and $\eta=1$): \\
In this case, it is enough to consider

\begin{equation*}
    \parallel u_0\parallel_{\dot{{\cal FB}}^{2-\frac{3}{p}}_{p,q}}<C_1\,\,\mbox{and}\,\,
    \parallel \theta_0\parallel_{\dot{{\cal FB}}^{-\frac{3}{p}}_{p,q}}<C_2.
\end{equation*}
\item [($a_2$)] ($\nu=1$ and $\eta<1$): \\
In this case, the initial condition is reduced to:

\begin{equation*}
    \parallel u_0\parallel_{\dot{{\cal FB}}^{2-\frac{3}{p}}_{p,q}}<C_1\cdot\eta\,\,\mbox{and}\,\,
    \parallel \theta_0\parallel_{\dot{{\cal FB}}^{-\frac{3}{p}}_{p,q}}<C_2\cdot \eta^{2}.
\end{equation*}
Observe that the norm $\parallel \theta_0\parallel_{\dot{{\cal FB}}^{-\frac{3}{p}}_{p,q}}$ could be taken smaller than the norm $\parallel u_0\parallel_{\dot{{\cal FB}}^{2-\frac{3}{p}}_{p,q}}$, if we consider $\eta$ close to zero.
\item [($a_3$)] ($\nu=1$ and $\eta>1$): \\
In this case, we can consider:

\begin{equation*}
    \parallel u_0\parallel_{\dot{{\cal FB}}^{2-\frac{3}{p}}_{p,q}}<C_1\,\,\mbox{and}\,\,
    \parallel \theta_0\parallel_{\dot{{\cal FB}}^{-\frac{3}{p}}_{p,q}}<C_2\cdot \eta.
\end{equation*}
Since $1<\eta$, it is enough to consider $\parallel \theta_0\parallel_{\dot{{\cal FB}}^{-\frac{3}{p}}_{p,q}}<C_2$. However, we can take large enough initial temperature disturbance $\parallel \theta_0\parallel_{\dot{{\cal FB}}^{-\frac{3}{p}}_{p,q}}$ if we consider large enough values for $\eta$.

\item [($b_1$)] ($\nu<1$ and $\eta=1$): \\
In this case, we can consider:

\begin{equation*}
    \parallel u_0\parallel_{\dot{{\cal FB}}^{2-\frac{3}{p}}_{p,q}}<C_1\cdot\nu^{2}\,\,\mbox{and}\,\,
    \parallel \theta_0\parallel_{\dot{{\cal FB}}^{-\frac{3}{p}}_{p,q}}<C_2\cdot \nu^{2}.
\end{equation*}
Observe that, if we consider small enough $\nu$, both the norms $\parallel u_0\parallel_{\dot{{\cal FB}}^{2-\frac{3}{p}}_{p,q}}$ and $\parallel \theta_0\parallel_{\dot{{\cal FB}}^{-\frac{3}{p}}_{p,q}}$
have the same quadratic smallness decaying conditions.

\item [($b_2$)] ($\nu<1$ and $\eta<1$): \\
In this case, it is enough to consider:

\begin{equation*}
    \parallel u_0\parallel_{\dot{{\cal FB}}^{2-\frac{3}{p}}_{p,q}}<C_1\cdot k^{2}\,\,\mbox{and}\,\,
    \parallel \theta_0\parallel_{\dot{{\cal FB}}^{-\frac{3}{p}}_{p,q}}<C_2\cdot k^{3},
\end{equation*}
where $k=\min(\nu,\eta)$. Thus, observe that as $k$ is close to zero then the norm for the temperature  $\parallel \theta_0\parallel_{\dot{{\cal FB}}^{-\frac{3}{p}}_{p,q}}$ could be taken smaller than the norm for the velocity  $\parallel u_0\parallel_{\dot{{\cal FB}}^{2-\frac{3}{p}}_{p,q}}$.

\item [($b_3$)] ($\nu<1$ and $\eta>1$): \\
In this case, it is enough to consider:

\begin{equation*}
    \parallel u_0\parallel_{\dot{{\cal FB}}^{2-\frac{3}{p}}_{p,q}}<C_1\cdot \nu^{2}\,\,\mbox{and}\,\,
    \parallel \theta_0\parallel_{\dot{{\cal FB}}^{-\frac{3}{p}}_{p,q}}<C_2\cdot\eta\cdot \nu^{2}.
\end{equation*}
Since $\nu<1<\eta$, it is enough to take 

\begin{equation*}
    \parallel u_0\parallel_{\dot{{\cal FB}}^{2-\frac{3}{p}}_{p,q}}<C_1\cdot \nu^{2}\,\,\mbox{and}\,\,
    \parallel \theta_0\parallel_{\dot{{\cal FB}}^{-\frac{3}{p}}_{p,q}}<C_2\cdot \nu^{3}.
\end{equation*}
 Thus, observe that as $\nu$ is close to zero then the norm for the temperature  $\parallel \theta_0\parallel_{\dot{{\cal FB}}^{-\frac{3}{p}}_{p,q}}$ could be taken smaller than the norm for the velocity  $\parallel u_0\parallel_{\dot{{\cal FB}}^{2-\frac{3}{p}}_{p,q}}$. Nevertheless, it is possible to consider another interesting subcases:
 \begin{itemize}
\item [($b_3.i$)] ($\eta=\frac{1}{\nu}$): In this subcase it is enough to have

\begin{equation*}
    \parallel u_0\parallel_{\dot{{\cal FB}}^{2-\frac{3}{p}}_{p,q}}<C_1\cdot \nu^{2}\,\,\mbox{and}\,\,
    \parallel \theta_0\parallel_{\dot{{\cal FB}}^{-\frac{3}{p}}_{p,q}}<C_2\cdot \nu.
\end{equation*}
Thus, as $\nu$ is taken close to zero, the value for the norm $\parallel u_0\parallel_{\dot{{\cal FB}}^{2-\frac{3}{p}}_{p,q}}$ could be taken smaller than the norm $\parallel \theta_0\parallel_{\dot{{\cal FB}}^{-\frac{3}{p}}_{p,q}}$.
\item [($b_3.ii$)] ($\eta=\frac{1}{\nu^{2}}$): In this subcase it is enough to have

\begin{equation*}
    \parallel u_0\parallel_{\dot{{\cal FB}}^{2-\frac{3}{p}}_{p,q}}<C_1\cdot \nu^{2}\,\,\mbox{and}\,\,
    \parallel \theta_0\parallel_{\dot{{\cal FB}}^{-\frac{3}{p}}_{p,q}}<C_2.
\end{equation*}
Observe that, as $\nu$ is small enough, the norm $\parallel u_0\parallel_{\dot{{\cal FB}}^{2-\frac{3}{p}}_{p,q}}$ should be taken small enough, however the smallness condition from the another norm  $\parallel \theta_0\parallel_{\dot{{\cal FB}}^{-\frac{3}{p}}_{p,q}}$ is independent of both the parameters $\nu$ and $\eta$. 

\item [($b_3.iii$)] ($\eta=\frac{1}{\nu^{3}}$): In this subcase it is enough to have

\begin{equation*}
    \parallel u_0\parallel_{\dot{{\cal FB}}^{2-\frac{3}{p}}_{p,q}}<C_1\cdot \nu^{2}\,\,\mbox{and}\,\,
    \parallel \theta_0\parallel_{\dot{{\cal FB}}^{-\frac{3}{p}}_{p,q}}<C_2\cdot\frac{1}{\nu}.
\end{equation*}
Observe that, if we consider small enough values for $\nu$, the norm for the initial velocity $\parallel u_0\parallel_{\dot{{\cal FB}}^{2-\frac{3}{p}}_{p,q}}$ should be taken small enough, however the norm for initial temperature  $\parallel \theta_0\parallel_{\dot{{\cal FB}}^{-\frac{3}{p}}_{p,q}}$ could be considered large enough.
 \end{itemize}
\item [($c_1$)] ($\nu>1$ and $\eta=1$): 
In this case we can consider:

\begin{equation*}
    \parallel u_0\parallel_{\dot{{\cal FB}}^{2-\frac{3}{p}}_{p,q}}<C_1\cdot \nu\,\,\mbox{and}\,\,
    \parallel \theta_0\parallel_{\dot{{\cal FB}}^{-\frac{3}{p}}_{p,q}}<C_2\cdot \nu.
\end{equation*}
Since $\nu>1$, it is enough to consider

 \begin{equation*}
    \parallel u_0\parallel_{\dot{{\cal FB}}^{2-\frac{3}{p}}_{p,q}}<C_1\,\,\mbox{and}\,\,
    \parallel \theta_0\parallel_{\dot{{\cal FB}}^{-\frac{3}{p}}_{p,q}}<C_2.
\end{equation*}
However, if we consider  large enough values for $\nu$ it is possible to take large enough values for both the norms $\parallel u_0\parallel_{\dot{{\cal FB}}^{2-\frac{3}{p}}_{p,q}}$ and $\parallel \theta_0\parallel_{\dot{{\cal FB}}^{-\frac{3}{p}}_{p,q}}$.
\item [($c_2$)] ($\nu>1$ and $\eta<1$): 
In this case we have:

\begin{equation*}
    \parallel u_0\parallel_{\dot{{\cal FB}}^{2-\frac{3}{p}}_{p,q}}<C_1\cdot \nu\cdot\eta\,\,\mbox{and}\,\,
    \parallel \theta_0\parallel_{\dot{{\cal FB}}^{-\frac{3}{p}}_{p,q}}<C_2\cdot \nu\cdot\eta^{2}.
\end{equation*}
Since $\eta<1<\nu$, it is enough to consider

 \begin{equation*}
    \parallel u_0\parallel_{\dot{{\cal FB}}^{2-\frac{3}{p}}_{p,q}}<C_1\cdot\eta^{2}\,\,\mbox{and}\,\,
    \parallel \theta_0\parallel_{\dot{{\cal FB}}^{-\frac{3}{p}}_{p,q}}<C_2\cdot\eta^{3}.
\end{equation*}
Observe that as $\eta$ is close to zero, it is possible to consider the value for the norm $\parallel \theta_0\parallel_{\dot{{\cal FB}}^{-\frac{3}{p}}_{p,q}}$ smaller than the value for the norm $\parallel u_0\parallel_{\dot{{\cal FB}}^{2-\frac{3}{p}}_{p,q}}$. Nevertheless, it is possible to consider another relevant subcases:
\begin{itemize}
\item [($c_2.i$)] ($\nu=\frac{1}{\eta}$): In this subcase it is enough to consider 

\begin{equation*}
    \parallel u_0\parallel_{\dot{{\cal FB}}^{2-\frac{3}{p}}_{p,q}}<C_1\,\,\mbox{and}\,\,
    \parallel \theta_0\parallel_{\dot{{\cal FB}}^{-\frac{3}{p}}_{p,q}}<C_2\cdot\eta.
\end{equation*}
This initial condition establish that for small enough values for $\eta$, the norm for the temperature $\parallel u_0\parallel_{\dot{{\cal FB}}^{2-\frac{3}{p}}_{p,q}}$ should be taken small enough, however, the smallness condition for the velocity does not depend of both the parameters $\nu$ and $\eta$.
\item [($c_2.ii$)] ($\nu=\frac{1}{\eta^{2}}$): In this subcase it is enough to consider 

\begin{equation*}
    \parallel u_0\parallel_{\dot{{\cal FB}}^{2-\frac{3}{p}}_{p,q}}<C_1\cdot\frac{1}{\eta}\,\,\mbox{and}\,\,
    \parallel \theta_0\parallel_{\dot{{\cal FB}}^{-\frac{3}{p}}_{p,q}}<C_2.
\end{equation*}
This initial condition establish that for small enough values for $\eta$, it is possible to consider large enough values for the norm $\parallel u_0\parallel_{\dot{{\cal FB}}^{2-\frac{3}{p}}_{p,q}}$ and the smallness condition for the temperature does not depends of both the parameters $\nu$ and $\eta$.

\item [($c_2.iii$)] ($\nu=\frac{1}{\eta^{3}}$): In this subcase it is enough to consider 

\begin{equation*}
    \parallel u_0\parallel_{\dot{{\cal FB}}^{2-\frac{3}{p}}_{p,q}}<C_1\cdot\frac{1}{\eta^{2}}\,\,\mbox{and}\,\,
    \parallel \theta_0\parallel_{\dot{{\cal FB}}^{-\frac{3}{p}}_{p,q}}<C_2\cdot\frac{1}{\eta}.
\end{equation*}
Observe that, as $\eta$ is close to zero, it is possible to consider both the norms $\parallel u_0\parallel_{\dot{{\cal FB}}^{2-\frac{3}{p}}_{p,q}}$ and  $\parallel \theta_0\parallel_{\dot{{\cal FB}}^{-\frac{3}{p}}_{p,q}}$ large enough. Particularly, the norm $\parallel u_0\parallel_{\dot{{\cal FB}}^{2-\frac{3}{p}}_{p,q}}$ could be taken bigger than the norm $\parallel \theta_0\parallel_{\dot{{\cal FB}}^{-\frac{3}{p}}_{p,q}}$.
\end{itemize}
Finally, we consider the last case
\item[($c_3$)]($\nu>1$ and $\eta>1$):
In this case it is enough to consider

\begin{equation*}
    \parallel u_0\parallel_{\dot{{\cal FB}}^{2-\frac{3}{p}}_{p,q}}<C_1\cdot \nu\,\,\mbox{and}\,\,
    \parallel \theta_0\parallel_{\dot{{\cal FB}}^{-\frac{3}{p}}_{p,q}}<C_2\cdot \nu\cdot\eta.
\end{equation*}
Since $1<\nu$ and $1<\eta$ it is enough to consider

\begin{equation*}
    \parallel u_0\parallel_{\dot{{\cal FB}}^{2-\frac{3}{p}}_{p,q}}<C_1\,\,\mbox{and}\,\,
    \parallel \theta_0\parallel_{\dot{{\cal FB}}^{-\frac{3}{p}}_{p,q}}<C_2.
\end{equation*}
Nevertheless, it is possible to consider large enough valued for the norms $\parallel u_0\parallel_{\dot{{\cal FB}}^{2-\frac{3}{p}}_{p,q}}$ and 
$ \parallel \theta_0\parallel_{\dot{{\cal FB}}^{-\frac{3}{p}}_{p,q}}$ if we consider large enough values for both the parameters $\nu$ and $\eta$.
\end{itemize}

\textbf{Acknowledgment}:\newline
The author was supported by Cnpq (Brazil) and by University of Campinas, SP, Brazil.

%\textbf{Acknowledgment}:\newline
%I would like to thanks Prof. Lucas C. F. Ferreira (UNICAMP) for fruitful %talks on this problem which certainly helped to establish the final form %of this text. This work was supported by University of Campinas (Brazil).
\end{section}

%%%%%%%%%%%%%%%%%%%%%%%%%%%%%%%%%%%%%%%%%%%%%%%%%%%%%%
%%%%%%%%%%%%%%%%%%%%%%%%%%%%%%%%%%%%%%%%%%%%%%%%%%%%%%%

%Leithold L. Aurazo-Alvarez\\
\emph{E-mail address}: leithold@unicamp.br


\begin{thebibliography}{99}
	
%%%%%%%%%%%
%%%%%%%%%%% A %%%%%%%%%

\bibitem{AbidiHmidiKeraaniGlobalregNSB}{\small H. Abidi, T. Hmidi and S. Keraani, On the global regularity of axisymmetric Navier-Stokes-Boussinesq system, Discrete and Continuous Dynamical Systems, 29, 3 (2011), 737--756.}
\bibitem{AbidiHmidiOntheGlobalBQ}{\small 
H. Abidi and T. Hmidi, On the global well-posedness for Boussinesq system. Journal of Differential Equations, 233, 1 (2007), 199--220.}


\bibitem {Almeida-Ferreira-Lima-UniformNSCFBMorrey}{\small
M.F. de Almeida, L.C.F. Ferreira and L.S.M. Lima, Uniform global well-posedness of the Navier-Stokes-Coriolis system in a new critical space. Math. Z., 287 (2017), 735--750.
}

\bibitem {Aurazo-Ferreira-Global-BC-Stratification-FBM-2021}{\small 
L.L. Aurazo-Alvarez and L.C.F. Ferreira,  Global well-posedness for the fractional Boussinesq-Coriolis system with stratification in a framework of Fourier-Besov type. Partial Differ. Equ. Appl., 2 (2021), 62}


%%%%%%%%%%% B %%%%%%%%%

\bibitem{BatchelorAnIntrodFD}{\small 
G. K. Batchelor, An introduction to fluid dynamics. Cambridge university press, 2000.}

\bibitem{GuoSpectralmethodNBsystem}{\small 
G. Boling, Spectral method for solving two-dimensional Newton-Boussinesq equations. Acta Mathematicae Applicatae Sinica, 5 (1989), 208--218.}	

\bibitem {BonyParaproduct}{\small
J.-M. Bony, Calcul symbolique et propagation des singularit\'es pour les \'equations aux d\'eriv\'ees partielles non lin\'eaires. Annales scientifiques de l'\'Ecole Normale Sup\'erieure, 4, 14, 2 (1981), 209--246.
}


%%%%%%%%%%% C %%%%%%%%%
\bibitem{Cannon-DiBenedetto}{\small 
J.R. Cannon and E. DiBenedetto, The initial value problem for the Boussinesq equations with data in Lp . In: Rautmann, R. (eds) Approximation Methods for Navier-Stokes Problems. Lecture Notes in Mathematics, Springer, Berlin, Heidelberg., (1980) 771.}


\bibitem{CaoWuGlobalRegVerticalDiss}{\small  C. Cao and J. Wu, Global Regularity for the Two-Dimensional Anisotropic Boussinesq Equations with Vertical Dissipation. Arch Rational Mech Anal, 208 (2013), 985--1004.}

\bibitem{ChaeGlobalRegBQPartViscosity}{\small D. Chae,
Global regularity for the 2D Boussinesq equations with partial viscosity terms. Advances in Mathematics, 203, 2 (2006), 497--513.}
%%%%%%%%%%% D %%%%%%%%%


\bibitem{DanchinPaicutheoremesLFK}{\small 
R. Danchin and M. Paicu, Les th\'eor\`emes de Leray et de Fujita-Kato pour le syst\`eme de Boussinesq partiellement visqueux. Bulletin de la Soci\'et\'e Math\'ematique de France, 136, 2 (2008), 261--309.}


\bibitem{DanchinPaicuExistUniqueLorentz}{\small R. Danchin and M. Paicu,
Existence and uniqueness results for the Boussinesq system with data in Lorentz spaces, Physica D: Nonlinear Phenomena, 237, 10-12 (2008), 1444--1460.}


\bibitem{DanchinPaicuGlobalIBQYudovichData}{\small R. Danchin and M. Paicu, Global Well-Posedness Issues for the Inviscid Boussinesq System with Yudovich's Type Data. Commun. Math. Phys., 290 (2009), 1--14.}



\bibitem{DanchinPaicuGlobalExAnisBQ}{\small R. Danchin and M. Paicu, Global existence results for the anisotropic Boussinesq system in dimension two. Mathematical Models and Methods in Applied Sciences, 21, 3 (2011), 421--457.}

%%%%%%%%%%% F %%%%%%%%%


\bibitem{FangLeZhangGlobalAxiNonzeroSwirl}{\small 
D. Fang, W. Le and T. Zhang, Global solutions of 3D axisymmetric Boussinesq equations with nonzero swirl. Nonlinear Analysis, 166 (2018), 48--86.}


\bibitem{Ferreira-Villamizar-R-Existence-convection-pseudomeasure-2008}{\small L.C.F. Ferreira, E.J. Villamizar-Roa, Existence of solutions to the convection problem in a pseudomeasure-type space. Proc. R. Soc. Lond. A Math. Phys. Eng. Sci., 464, 2096 (2008), 1983--1999.
}

\bibitem{Ferreira-Villamizar-R-Stability-Boussinesq-Weak-Lp-2010}{\small L.C.F. Ferreira, E.J. Villamizar-Roa, On the stability problem for the Boussinesq equations in weak-$L^p$ spaces, Communications on Pure and Applied Analysis, 9, 3 (2010), 667--684.
}

\bibitem{Ferreira-Almeida-Wellposedness-Boussinesq-Morrey-2011}
{\small L.C.F. Ferreira and M.F. de Almeida, On the well posedness and large-time behavior for Boussinesq equations in Morrey spaces. Differential Integral Equations, 24, 7-8 (2011), 719--742.
}


\bibitem{Ferreira-Perez-Villamizar-B-Lorentz-M-SBQ-2018}
{\small L.C.F. Ferreira, J.E. P\'erez-L\'opez and E.J. Villamizar-Roa, On the product in Besov-Lorentz-Morrey spaces and existence of solutions for the stationary Boussinesq equations. Communications on Pure and Applied Analysis, 17, 6 (2018), 2423--2439. 
}

\bibitem{Ferreira-Precioso-Existence-micropolar-Besov-Morrey-2013}
{\small L.C.F. Ferreira and J.C. Precioso, Existence of solutions for the 3D-micropolar fluid system with initial data in Besov-Morrey spaces. Z. Angew. Math. Phys., 64 (2013), 1699--1710.
}



%%%%%%%%%%% H %%%%%%%%%


\bibitem{HmidiKeraaniOntheGlobalBQZeroDiffus}{\small T. Hmidi and S. Keraani, On the global well-posedness of the two-dimensional Boussinesq system with a zero diffusivity. Adv. Differential Equations, 12, 4 (2007), 461--480.}

\bibitem{HmidiRoussetGlobalwellNSBQaxi}{\small T. Hmidi and F. Rousset, Global well-posedness for the Navier-Stokes-Boussinesq system with axisymmetric data. Annales de l'I.H.P. Analyse non lin\'eaire, 27, 5 (2010), 1227--1246.}

\bibitem{Hou-LiGlobalWVBQ}{\small 
T. Y. Hou and  C. Li, Global well-posedness of the viscous Boussinesq equations, 12, 1 (2005), 1--12.}
%%%%%%%%%%% K %%%%%%%%%

\bibitem {Konieczny-Yoneda2011}{\small P. Konieczny and T. Yoneda, On dispersive effect of the Coriolis force for the stationary Navier-Stokes equations. Journal of Differential Equations, 250, 10 (2011), 3859--3873.}
%%%%%%%%%%% L %%%%%%%%%

\bibitem {Lemarie2002}{\small P. G. Lemarie-Rieusset, Recent developments in the Navier-Stokes equations, Chapman and Hall, Research Notes in Maths., 431, 2002.}

\bibitem{LiTitiGlobalBQVerticalDiss}{\small  J. Li and E.S. Titi, Global Well-Posedness of the 2D Boussinesq Equations with Vertical Dissipation. Arch Rational Mech Anal, 220 (2016), 983--1001.}


%%%%%%%%%%% M %%%%%%%%%

\bibitem{MajdaIntrodPDEWAO}{\small 
A. Majda, Introduction to PDEs and Waves for the Atmosphere and Ocean. American Mathematical Soc., 2003.}

\bibitem{MiaoZhengBQHorizontalDissipation}{\small 
C. Miao and X. Zheng, On the Global Well-posedness for the Boussinesq System with Horizontal Dissipation. Commun. Math. Phys., 321 (2013), 33--67.}

\bibitem{MiaoZhengBQAxihorizontalViscosity}{\small C. Miao and X. Zheng,
Global well-posedness for axisymmetric Boussinesq system with horizontal viscosity,
Journal de Math\'ematiques Pures et Appliqu\'es, 101, 6 (2014), 842--872.}


%%%%%%%%%%% P %%%%%%%%%

\bibitem{PedloskyGeophysicalFD}{\small 
J. Pedlosky, Geophysical Fluid Dynamics. Springer New York, NY, 1987.}

\end{thebibliography}
\end{document}